\documentclass[smallextended,envcountsect]{svjour3}
\smartqed
\usepackage{algorithm}
\usepackage{graphicx}

\usepackage{amsmath,amsfonts,bm}

\def\1{\bm{1}}

\def\ve{{\bm{e}}}

\def\vg{{\bm{g}}}

\def\vu{{\bm{u}}}
\def\vv{{\bm{v}}}

\def\vx{{\bm{x}}}
\def\vy{{\bm{y}}}

\DeclareMathAlphabet{\mathsfit}{\encodingdefault}{\sfdefault}{m}{sl}
\SetMathAlphabet{\mathsfit}{bold}{\encodingdefault}{\sfdefault}{bx}{n}

\newcommand{\R}{\mathbb{R}}

\DeclareMathOperator{\sign}{sign}

\usepackage[hidelinks]{hyperref}

\usepackage{url}
\usepackage{booktabs}
\usepackage{algorithmic} 
\usepackage{tikz}
\usetikzlibrary{arrows.meta,calc}

\spnewtheorem{assumption}[theorem]{Assumption}{\bfseries}{\rmfamily}

\providecommand{\sign}{\operatorname{sign}}
\providecommand{\Sign}{\operatorname{Sign}}

\newcommand{\assref}[1]{Assumption~\ref{#1}}
\newenvironment{assumptionbox}[2][]{%
  \begin{assumption}[#1]\label{#2}}{\end{assumption}}

\usepackage{tabularx,booktabs}  

\title{Norm-Constrained Flows and Sign-Based Optimization: Theory and Algorithms}

\author{Valentin Leplat \and Sergio Mayorga \and Roland Hildebrand \and Alexander Gasnikov}

\institute{Valentin Leplat, Corresponding author \at
Department of Mathematics and Operational Research, Universit\'e de Mons, Rue de Houdain 9, 7000 Mons, Belgium\\
\email{valentin.leplat@gmail.com}
\and
Sergio Mayorga \at
Institute of Data Science and Artificial Intelligence, Innopolis University, Innopolis, Russia\\
\email{s.mayorga@innopolis.ru}
\and
Roland Hildebrand \at
Department of Mathematical Foundations of Control, Moscow Institute of Physics and Technology, Moscow, Russia\\
\email{khildebrand.r@mipt.ru}
\and
Alexander Gasnikov \at
Department of Mathematical Foundations of Control, Moscow Institute of Physics and Technology, Moscow, Russia\\
\email{gasnikov@yandex.ru}
}

\date{}

\begin{document}

\maketitle

\vspace{-2.0cm}
\begin{abstract}
Sign Gradient Descent (SignGD) uses only the coordinate-wise sign of the gradient. We study this method through norm-constrained continuous-time dynamics: at each point, the velocity is chosen to minimize the directional derivative over a unit norm ball. This recovers the sign flow for the $\ell_\infty$ constraint, normalized gradient flow for $\ell_2$, and greedy coordinate directions for $\ell_1$. We formulate the resulting dynamics as a set-valued differential inclusion, prove existence of solutions, and derive an exact energy identity. This identity also gives finite-time convergence under a Polyak--\L{}ojasiewicz inequality expressed in the dual norm. We then connect the canonical set-valued flow with classical Filippov regularization of discontinuous selectors, which gives a precise description of crossing and sliding near switching sets. Motivated by this behavior, we introduce two face-aware SignGD variants, one-hit freeze and two-hit sliding-track. Both methods modify the usual sign direction by damping selected coordinates when a crossing or persistent switching is detected. We derive descent certificates for these damped sign directions and introduce a safeguard that preserves a global linear convergence rate under $\ell_\infty$-smoothness and a Polyak--\L{}ojasiewicz inequality. For the $\ell_1$ flow, we also introduce convex-combination updates on active faces and prove a corresponding linear convergence guarantee. Finally, we analyze mass-aware face restrictions and an inertial SignGD method with restart. Numerical experiments illustrate the certified step rules, the face-aware variants, and the inertial scheme.
\end{abstract}
\keywords{Sign gradient descent \and norm-constrained flows \and Filippov differential inclusions \and sliding modes \and strongly convex optimization}

\section{Introduction}

Sign-based optimization methods use simple updates and are inexpensive to communicate. They are also insensitive to the magnitude of individual gradient coordinates. The basic deterministic SignGD iteration is
\begin{equation}
    \vx_{k+1}=\vx_k-\eta_k\,\operatorname{sign}\!\bigl(\nabla f(\vx_k)\bigr),
\end{equation}
where $\eta_k>0$ and the sign is taken componentwise. Sign-based methods first became prominent in stochastic and distributed optimization, where one-bit communication and robustness to noisy gradients are important, see~\cite{bernstein2018signsgdcompressedoptimisationnonconvex} for more details.

There is also a deterministic geometric interpretation. SignGD is steepest descent when velocity is measured in the $\ell_\infty$ norm. This viewpoint was developed by Balles, Pedregosa, and Le Roux (arXiv 2002.08056, 2020), and linear convergence of scaled SignGD was proved by Li et al.~\cite{li2021faster} under strong convexity and PL-type assumptions. We build on this foundation. Our focus is the set-valued geometry that appears when the norm-steepest direction is not unique, the Filippov dynamics generated by discontinuous selectors, and the discrete algorithms suggested by this geometry.

The main contributions are as follows.
\begin{itemize}
  \item \textbf{Norm-constrained flows.} We define the optimal velocity set
  $\mathcal V(x)=\arg\min_{\|v\|\le1}\langle\nabla f(x),v\rangle=-\partial\|\cdot\|_*(\nabla f(x))$. This single formula gives the $\ell_\infty$, $\ell_2$, and $\ell_1$ flows, including zeros and ties. We prove existence of trajectories and an exact energy identity along every admissible trajectory.
  \item \textbf{Continuous-time convergence.} The energy identity gives convergence estimates for convex objectives and finite-time convergence under a norm-PL inequality. For the $\ell_\infty$ sign flow, Ma et al.~\cite{ma2022dynamic} established finite-time convergence under a PL inequality. Our statement covers arbitrary norms and set-valued face dynamics.
  \item \textbf{Filippov switching and sliding.}
We relate the canonical optimization inclusion to the classical Filippov
regularization of discontinuous selectors. Away from switching sets, the
Filippov dynamics coincide with the ordinary selector dynamics. On switching
sets, the admissible velocities are convex combinations of locally attainable
norm-steepest directions. On regular switching hypersurfaces, we give a simple
criterion that distinguishes crossing from sliding.

\item \textbf{Face-aware SignGD.}
Motivated by the Filippov dynamics, we introduce two modifications of SignGD:
the one-hit freeze rule and the two-hit sliding-track rule. The one-hit rule
sets the velocity of a coordinate to zero immediately after a genuine sign
crossing. The two-hit rule detects two consecutive crossings and replaces the
usual velocity $\pm1$ by a coefficient in $[-1,1]$ that approximates a sliding
velocity. Both methods therefore generate coordinatewise damped sign
directions of the form
\[
d(w,g)=w\odot\operatorname{sign}(g),
\qquad 0\le w_i\le1.
\]
For the $\ell_1$ flow, we also introduce updates obtained by taking convex
combinations of the active coordinate directions.

\item \textbf{Convergence guarantees for the face-aware updates.}
We analyze the damped sign directions generated by the one-hit and two-hit
rules. We first derive a general descent bound and then introduce a safeguard
that accepts a face-aware direction whenever it preserves enough descent and
otherwise uses the ordinary scaled SignGD direction. Under
$\ell_\infty$-smoothness and a PL inequality, the safeguarded method converges
linearly. For the $\ell_1$ flow, we prove a corresponding linear rate for every
convex-combination update under the appropriate norm-smoothness and PL
assumptions.
  \item \textbf{SignGD, active faces, and inertia.} We give a self-contained deterministic analysis of scaled SignGD under separable smoothness and strong convexity. We then refine the rate by accounting for the retained gradient mass on a face. Finally, we analyze an inertial SignGD method with restart; it keeps the global linear guarantee of the basic adaptive method and quantifies the extra one-step gain obtained by successful extrapolation.
\end{itemize}

We conclude with numerical experiments comparing the basic, face-aware, and inertial methods on smooth strongly convex problems. For the adaptive rules, the theoretical constants must come from a genuine separable upper model (or an equivalent $\ell_\infty$ smoothness bound); diagonal Hessian entries alone are not sufficient for a generic dense objective.

\subsection{Notation and Assumptions}
We write $\vx\in\R^d$ for parameters and $f:\R^d\to\R$ for the objective.
Norms $\|\cdot\|_p$ are standard; $\|\cdot\|_\ast$ denotes the dual norm.

For $\vu\in\R^d$, the \emph{element-wise} sign used in discrete algorithms is
\[
  \sign(\vu)_i=\begin{cases}
  1 & u_i>0,\\[2pt]
  0 & u_i=0,\\[2pt]
  -1 & u_i<0,
  \end{cases}
  \qquad i=1,\ldots,d.
\]
For continuous-time arguments we use the \emph{set-valued} sign
$\Sign(\vu)\subset[-1,1]^d$ defined componentwise by
\[
  \Sign(\vu)_i=\begin{cases}
  \{1\} & u_i>0,\\[-2pt]
  [-1,1] & u_i=0,\\
  \{-1\} & u_i<0.
  \end{cases}
\]
Thus a zero coordinate does not determine a unique velocity. This is exactly what makes the continuous-time model set-valued. Section~\ref{sec:filippov} explains how this inclusion is related to classical Filippov regularization of discontinuous single-valued selectors.
We denote by $\ve^{(i)}$ the $i$-th standard basis vector and by
$\partial_i f(\vx)$ the $i$-th component of $\nabla f(\vx)$.

\begin{assumptionbox}[Standing assumptions]{ass:standing}

The function $f$ is continuously differentiable.
\end{assumptionbox}

\begin{assumptionbox}[Separable smoothness]{ass:coord}
There are constants $L_i>0$, $i=1,\ldots,d$, such that for all $x,y$,
\[
f(y) \le f(x) + \langle \nabla f(x), y-x\rangle
       + \tfrac12 \sum_{i=1}^d L_i (y_i-x_i)^2.
\]
\end{assumptionbox}

\noindent\emph{Note.} Assumption~\ref{ass:coord} is the separable quadratic upper model used in several analyses of sign methods; see, for example, Balles, Pedregosa, and Le Roux (arXiv 2002.08056, 2020). It is stronger than merely assuming that each partial derivative is Lipschitz when the other coordinates are fixed. It implies the Euclidean quadratic upper model with constant $L=\max_i L_i$ and the $\ell_\infty$ quadratic upper model with constant $\sum_iL_i$. We use only these upper bounds below; no stronger claim about gradient Lipschitz continuity is needed.

\begin{assumptionbox}[Strong convexity (when stated)]{ass:strong}
$f$ is $\mu$--strongly convex for some $\mu>0$, i.e.,
\[
f(\vy)\ \ge\ f(\vx) + \langle \nabla f(\vx), \vy-\vx\rangle + \tfrac{\mu}{2}\|\vy-\vx\|_2^2
\quad\text{for all }\vx,\vy\in\R^d.
\]
\end{assumptionbox}

We next review related work on sign-based methods, continuous-time dynamics, and non-Euclidean optimization.


\section{Related Work}

\paragraph{Sign descent and non-Euclidean geometry.}
The communication benefits of SignSGD were emphasized by Bernstein et al.~\cite{bernstein2018signsgdcompressedoptimisationnonconvex}, followed by a large literature on error feedback, quantization, and stochastic sign methods~\cite{pmlr-v97-karimireddy19a,stich2018sparsifiedsgdmemory,alistarh2017qsgdcommunicationefficientsgdgradient}. For the deterministic method, Balles, Pedregosa, and Le Roux (arXiv 2002.08056, 2020) study SignGD as steepest descent in the $\ell_\infty$ norm and distinguish the separable quadratic upper model from the weaker $\ell_\infty$-smoothness assumption. Li et al.~\cite{li2021faster} analyze scaled SignGD and obtain linear convergence for strongly convex and PL objectives. These results provide the baseline deterministic theory for SignGD. We use the same geometry to study switching, sliding, active faces, and inertial variants.

\paragraph{Continuous-time sign dynamics.}
Continuous-time models have long been useful for understanding optimization algorithms and momentum~\cite{su2014differential,wibisono2016variational}. Ma et al.~\cite{ma2022dynamic} show that a continuous-time limit of RMSprop and Adam approaches the SignGD flow and prove finite-time convergence of that flow under a PL condition. Here we start from the norm-constrained optimal-velocity problem. This gives a set-valued inclusion for any norm, together with an exact descent identity that holds for every admissible velocity, including when the optimal direction is not unique.

\paragraph{Filippov systems.}
The standard theory of differential equations with discontinuous right-hand sides is developed in~\cite{Filippov1988}; a concise optimization-friendly tutorial is~\cite{Cortes2008Discontinuous}. The classical Filippov regularization removes arbitrary measure-zero subsets before taking the local closed convex hull. This detail matters at degenerate switching points. We therefore distinguish the canonical optimization inclusion from the Filippov regularization of a particular discontinuous selector. The two coincide only when all extreme velocities of the canonical set are locally attainable by that selector.

\paragraph{Momentum and stochastic analyses.}
Momentum can substantially change the behavior of normalized and sign methods. Cutkosky and Mehta~\cite{pmlr-v119-cutkosky20b} analyze normalized SGD with momentum, and Sun et al.~\cite{pmlr-v202-sun23l} give convergence results for SignSGD with momentum under weaker assumptions. Recent preprints also study sign operators under heavy-tailed noise (Kornilov et al., arXiv 2502.07923, 2025; Yu et al., arXiv 2602.07425, 2026) and non-Euclidean stochastic-gradient methods that include SignSGD, Lion, and Muon as particular cases (Kovalev and Borodich, arXiv 2511.11466, 2025). These stochastic results are complementary to the deterministic switching and sliding questions considered here.

\paragraph{Norm-ball LMOs, conditional gradients, and Muon.}
Pethick et al.~\cite{pethick2025lmo} study stochastic optimization methods based on linear minimization oracles over norm balls. Their framework includes Muon-type spectral-norm updates and shows how weight decay connects these unconstrained updates with stochastic conditional-gradient methods. Their arXiv version appeared in February 2025, before the non-Euclidean trust-region analysis of gradient orthogonalization by Kovalev (arXiv 2503.12645, 2025). This work therefore already makes explicit the connection between norm-ball LMOs, Muon, and conditional-gradient methods. Our focus is different: we study velocity-constrained differential inclusions, the geometry of their switching sets, and discrete updates motivated by crossing and sliding.

\paragraph{Positioning.}
The $\ell_\infty$ steepest-descent interpretation and the linear convergence theory of scaled SignGD provide an established baseline. Our focus is the full set-valued flow, its relation with Filippov regularization of discontinuous selectors, the distinction between crossing and sliding, and face-aware updates with convergence guarantees.

\section{Gradient Flow Perspective and Geometry}
\label{sec:flow-geometry}

\subsection{From Steepest Descent to an $\ell_\infty$-Constrained Flow}
Let $f:\R^d\to\R$ be continuously differentiable. The classical gradient flow minimizes $f(\vx(t))$ via
\begin{equation}
    \dot{\vx}(t) = -\nabla f(\vx(t)).
    \label{eq:unconstrained-flow}
\end{equation}
It arises from the steepest descent principle in $\ell_2$ geometry:
\begin{equation}
    \vv^\star \!=\! \arg\min_{\vv\in\R^d}\Big\{ \langle \nabla f(\vx(t)), \vv\rangle + \tfrac{1}{2}\|\vv\|_2^2 \Big\}
    \quad\Rightarrow\quad \vv^\star = -\nabla f(\vx(t)).
    \label{eq:steepest-descent}
\end{equation}

To impose a trust region, we replace the penalty with a hard norm constraint. In $\ell_\infty$ geometry we solve
\begin{equation}
    \min_{\|\vv\|_\infty \le 1}\ \langle \nabla f(\vx(t)), \vv\rangle.
    \label{eq:linf-problem}
\end{equation}
Write $p=\nabla f(\vx(t))$. The problem separates over the coordinates. If $p_i\neq0$, every minimizer satisfies
$v_i=-\operatorname{sign}(p_i)$. If $p_i=0$, any $v_i\in[-1,1]$ is optimal. Hence the full optimal-velocity set is
$-\Sign(p)$, and the corresponding \emph{sign gradient flow} is
\begin{equation}
    \dot{\vx}(t) \in -\operatorname{Sign}\!\big(\nabla f(\vx(t))\big).
    \label{eq:sign-flow}
\end{equation}
This is a differential inclusion whenever at least one partial derivative is zero.

A single forward-Euler step then gives the discrete algorithm, with the stepsize acting as an $\ell_\infty$ trust-region radius.

\subsection{Discretization: Sign Gradient Descent}
Applying one forward-Euler step with step size $\eta>0$ to~\eqref{eq:sign-flow} gives
\begin{equation}
    \vx_{k+1} = \vx_k - \eta\, \operatorname{sign}\!\big(\nabla f(\vx_k)\big),
\end{equation}
which is exactly the SignGD update. In this discretization, $\eta$ plays a dual role: it is both the time step and the effective $\ell_\infty$ trust-region radius $\|\vx_{k+1}-\vx_k\|_\infty\le \eta$.

\subsection{Alternative $\ell_1$-Constrained Flow and Sparse Updates}
Consider instead
\begin{equation}
    \dot{\vx}(t) \in \arg\min_{\|\vv\|_1 \le 1}\ \langle \nabla f(\vx(t)), \vv\rangle.
    \label{eq:l1-problem}
\end{equation}
Let $g=\nabla f(\vx(t))$ and define
\[
\mathcal I(\vx(t)):=\arg\max_i |g_i|.
\]
At a stationary point this set contains all coordinates.
The minimum value is $-\|g\|_\infty$. If $g\neq0$, the complete set of minimizers is
\begin{equation}
\mathcal V_1(\vx)
=\operatorname{conv}\Bigl\{-\operatorname{sign}(\partial_i f(\vx))\,\ve^{(i)}:
 i\in\mathcal I(\vx)\Bigr\}.
\label{eq:l1-explicit-V}
\end{equation}
If $\nabla f(\vx)=0$, every velocity in the $\ell_1$ unit ball is optimal:
\begin{equation}
\mathcal V_1(\vx)=B_1:=\{v\in\mathbb R^d:\|v\|_1\le1\}.
\label{eq:l1-zero-V}
\end{equation}
This makes the sparsity explicit. When the largest gradient magnitude is attained at a unique coordinate $i$, the optimal velocity has only one nonzero component,
\[
\dot{\vx}(t)=-\operatorname{sign}(\partial_i f(\vx(t)))\,\ve^{(i)}.
\]
A forward-Euler step is therefore the usual greedy coordinate update
\[
\vx_{k+1}=\vx_k-\eta\,\operatorname{sign}(\partial_i f(\vx_k))\,\ve^{(i)}.
\]
If several coordinates attain the same largest absolute gradient value, they \emph{tie}; any convex combination of their signed coordinate directions is optimal. This is the set-valued case studied below.

\subsection{A Unified View via Dual Norms}
\label{subsec:dual-norm-lens}

\begin{lemma}[Steepest descent under a norm constraint]\label{lemma:unified}
  Let $\|\cdot\|$ be any norm with dual $\|\cdot\|_\ast$, and let $g\in\R^d$.
Then
\[
\min_{\|\vv\|\le 1}\ \langle g,\vv\rangle \;=\; -\|g\|_\ast,
\qquad
\arg\min_{\|\vv\|\le 1}\ \langle g,\vv\rangle \;=\; -\,\partial\|\cdot\|_\ast(g),
\]
so the constrained steepest-descent flow can be written as
\[
\dot{\vx}(t)\in -\,\partial\|\cdot\|_\ast\!\big(\nabla f(\vx(t))\big).
\]  
\end{lemma}

\begin{proof}
    Deferred to Appendix~\ref{app:dual-norm-lemma} (Lemma A.1).
\end{proof}

\paragraph{Motivation.}
We choose a velocity $\vv$ inside the unit ball $\{\|\vv\|\le 1\}$ to decrease $f$
as fast as possible locally, i.e., to minimize the directional derivative
$\langle \nabla f,\vv\rangle$. The dual norm
$\|g\|_\ast=\sup_{\|\vv\|\le 1}\langle g,\vv\rangle$ measures the largest increase a linear form $g$ can induce on that ball; hence the largest
decrease is negative:
\[
\min_{\|\vv\|\le 1}\langle \nabla f,\vv\rangle
= -\,\|\nabla f\|_\ast.
\]
The minimizers form the face of the unit ball on which
$\langle \nabla f,\vv\rangle$ is smallest. Equivalently,
$\vv^\star\in-\,\partial\|\cdot\|_\ast(\nabla f)$.
When this face is a single point, the direction is unique. At zeros or ties,
the face can contain several directions, which makes the dynamics set-valued.

\paragraph{Discussion (concrete instances).}
Lemma~\ref{lemma:unified} recovers three standard algorithms from the same principle:
\begin{itemize}
  \item \textbf{$\ell_2$ geometry (Normalized GD).} Dual is $\ell_2$.
  For $g\neq 0$, $\partial\|\cdot\|_2(g)=\{g/\|g\|_2\}$, so
  $\dot{\vx}=-\nabla f/\|\nabla f\|_2$ and Euler discretization gives
  $\vx_{k+1}=\vx_k-\eta\,\nabla f(\vx_k)/\|\nabla f(\vx_k)\|_2$. At a stationary point the canonical velocity set is the $\ell_2$ unit ball; the usual discrete method simply stops.

  \item \textbf{$\ell_\infty$ geometry (Sign flow / SignGD).} Dual is $\ell_1$.
  $\partial\|\cdot\|_1(g)$ is the element-wise sign (with $[-1,1]$ at zeros),
  hence $\dot{\vx}\in-\,\operatorname{Sign}(\nabla f)$ and Euler discretization
  is $\vx_{k+1}=\vx_k-\eta\,\operatorname{sign}(\nabla f(\vx_k))$.

  \item \textbf{$\ell_1$ geometry (Greedy coordinate descent).} Dual is $\ell_\infty$.
  If $g\neq0$ and $I(g)=\arg\max_i|g_i|$, then
  \[
  \partial\|\cdot\|_\infty(g)
  =\operatorname{conv}\{\operatorname{sign}(g_i)\,\ve^{(i)}:i\in I(g)\}.
  \]
  Each extreme direction has exactly one nonzero component. If $|I(g)|>1$, several coordinates tie for the largest gradient magnitude and every convex combination of these extreme directions is optimal. At $g=0$, the subdifferential is the full $\ell_1$ unit ball.
\end{itemize}

\noindent

\textit{Trust-region view.}
With a discrete step size $\eta$, the subproblem becomes
$\min_{\|\vv\|\le \eta}\langle \nabla f(\vx_k),\vv\rangle$, and its minimizer set is
$-\eta\,\partial\|\cdot\|_\ast(\nabla f(\vx_k))$. Thus $\eta$ is both the Euler time step and the trust-region radius in the chosen norm.

Table~\ref{tab:ODE} summarizes the discussion above.

\begin{table}[ht!]
\caption{Optimization methods with constraints (the $\ell_2$ row uses the normalized flow)}
\label{tab:ODE}
\centering
\small
\resizebox{\textwidth}{!}{%
\begin{tabular}{llll}
\toprule
\textbf{Constraint} & \textbf{Flow} & \textbf{Euler update} & \textbf{Method} \\
\midrule
None & $\dot{\vx}=-\nabla f(\vx)$ & $\vx_{k+1}=\vx_k-\eta \nabla f(\vx_k)$ & Gradient Descent \\
$\ell_2$: $\|\vv\|_2\le 1$ &
$\dot{\vx}=- \frac{\nabla f(\vx)}{\|\nabla f(\vx)\|_2}$ &
$\vx_{k+1}=\vx_k-\eta\, \frac{\nabla f(\vx_k)}{\|\nabla f(\vx_k)\|_2}$ &
Normalized GD \\
$\ell_\infty$: $\|\vv\|_\infty\le 1$ &
$\dot{\vx}\in-\Sign(\nabla f(\vx))$ &
$\vx_{k+1}=\vx_k-\eta\, \sign(\nabla f(\vx_k))$ &
SignGD \\

$\ell_1$: $\|\vv\|_1\le 1$ & $\dot{\vx}\in \mathcal V_1(\vx)$ & $\vx_{k+1}=\vx_k-\eta\,\sign(\partial_i f(\vx_k))\,\ve^{(i)}$ & Greedy CD \\
\multicolumn{4}{l}{\small where $i\in\arg\max_j|\partial_j f(\vx_k)|$; see~\eqref{eq:l1-explicit-V}--\eqref{eq:l1-zero-V} for the full set-valued flow.}\\
\bottomrule
\end{tabular}%
}
\end{table}
\noindent\emph{Stationary points.} The normalized formula in the $\ell_2$ row is written for $\nabla f\neq0$. In the discrete methods we stop when $\nabla f=0$; the canonical set-valued formulation below also covers that case.

\subsection{Canonical Set-Valued Flow and an Exact Descent Identity}
\label{subsec:canonical-flow}
The minimizer in Lemma~\ref{lemma:unified} need not be unique. It is therefore natural to work with the whole optimal face rather than choosing a selector. For a norm $\|\cdot\|$ with dual norm $\|\cdot\|_*$, define
\begin{equation}
\mathcal V(x):=\operatorname*{argmin}_{\|v\|\le1}\langle\nabla f(x),v\rangle
=-\partial\|\cdot\|_*\bigl(\nabla f(x)\bigr).
\label{eq:canonical-V}
\end{equation}
We call
\begin{equation}
\dot x(t)\in\mathcal V(x(t))
\label{eq:canonical-flow}
\end{equation}
the \emph{canonical norm-constrained flow}. For $\ell_\infty$, this is $\dot x\in-\Sign(\nabla f(x))$. For $\ell_1$, it is
$\dot x\in-\partial\|\cdot\|_\infty(\nabla f(x))$, exactly the set $\mathcal V_1$ described in~\eqref{eq:l1-explicit-V}--\eqref{eq:l1-zero-V}.

\begin{proposition}[Existence of canonical norm-constrained trajectories]
\label{prop:canonical-existence}
Under Assumption~\ref{ass:standing}, the map $\mathcal V$ has nonempty, convex, compact values, is outer semicontinuous, and is globally bounded. Hence, for every initial point $x_0$, the inclusion~\eqref{eq:canonical-flow} admits at least one absolutely continuous solution starting from $x_0$.
\end{proposition}
\begin{proof}
For every $g$, $-\partial\|\cdot\|_*(g)$ is a nonempty compact convex subset of the primal unit ball. The graph of the subdifferential of a finite convex function is closed. Since $\nabla f$ is continuous, the graph of $\mathcal V$ is closed, and compact-valued local boundedness gives outer semicontinuity. The values are contained in the unit ball, so the map is globally bounded. Standard existence results for differential inclusions with such Marchaud maps apply; see~\cite{AubinCellina1984}.
\end{proof}

The optimal-velocity formulation also gives an exact identity along every trajectory, including on switching faces.
\begin{theorem}[Exact descent along the norm-constrained flow]
\label{thm:flow-energy}
Let $x(\cdot)$ be any absolutely continuous solution of~\eqref{eq:canonical-flow}. Then
\begin{equation}
\frac{d}{dt}f(x(t))=-\|\nabla f(x(t))\|_*
\qquad\text{for a.e. }t.
\label{eq:flow-energy}
\end{equation}
\end{theorem}
\begin{proof}
The chain rule for absolutely continuous curves gives, for a.e. $t$,
\[
\frac{d}{dt}f(x(t))
=
\langle\nabla f(x(t)),\dot x(t)\rangle.
\]
Every $v\in\mathcal V(x)$ minimizes the linear form over the unit ball, so
\[
\langle\nabla f(x),v\rangle=-\|\nabla f(x)\|_*.
\]
The same identity holds when $\nabla f(x)=0$.
\end{proof}

\subsection{Two Consequences for Convex and PL Objectives}
\label{subsec:flow-rates}
The first estimate only uses convexity. Assume that $f$ has a minimizer $x^\star$, write $f^\star:=f(x^\star)$, and suppose that the initial sublevel set satisfies
\begin{equation}
R_0:=\sup\{\|x-x^\star\|: f(x)\le f(x(0))\}<+\infty.
\label{eq:R0}
\end{equation}
Since~\eqref{eq:flow-energy} gives
\[
\frac{d}{dt}f(x(t))
=
-\|\nabla f(x(t))\|_*
\le 0
\qquad\text{for a.e. }t,
\]
the function \(t\mapsto f(x(t))\) is nonincreasing. Hence the trajectory
remains in the initial sublevel set,
\[
f(x(t))\le f(x(0))
\qquad\text{for all }t\ge 0.
\]
By the definition of \(R_0\), this implies
\[
\|x(t)-x^\star\|\le R_0.
\]

Using convexity of \(f\), we have
\[
f(x(t))-f^\star
\le
\langle \nabla f(x(t)),x(t)-x^\star\rangle.
\]
Applying the dual-norm inequality then gives
\[
f(x(t))-f^\star
\le
\|\nabla f(x(t))\|_*
\|x(t)-x^\star\|
\le
R_0\|\nabla f(x(t))\|_*.
\]
Therefore,
\[
\|\nabla f(x(t))\|_*
\ge
\frac{f(x(t))-f^\star}{R_0}.
\]

Let
\[
\Delta(t):=f(x(t))-f^\star.
\]
Combining the previous estimate with~\eqref{eq:flow-energy}, we obtain
\[
\dot{\Delta}(t)
=
-\|\nabla f(x(t))\|_*
\le
-\frac{1}{R_0}\Delta(t)
\qquad\text{for a.e. }t.
\]
Gr\"onwall's inequality then yields
\begin{equation}
\Delta(t)
\le
e^{-t/R_0}\Delta(0),
\end{equation}
or, equivalently,
\begin{equation}
f(x(t))-f^\star
\le
e^{-t/R_0}
\bigl(f(x(0))-f^\star\bigr).
\label{eq:convex-flow-exp}
\end{equation}

This bound mainly illustrates the behavior of the norm-constrained flow on a bounded sublevel set.

Now suppose that \(f\) satisfies the norm-PL inequality
\begin{equation}
\frac{1}{2}\|\nabla f(x)\|_*^2
\ge
\mu_*\bigl(f(x)-f^\star\bigr)
\qquad\text{for all }x,
\label{eq:norm-PL}
\end{equation}
for some \(\mu_*>0\). Let
\[
\Delta(t):=f(x(t))-f^\star.
\]
By Theorem~\ref{thm:flow-energy},
\[
\dot{\Delta}(t)
=
-\|\nabla f(x(t))\|_*
\qquad\text{for a.e. }t.
\]
Using~\eqref{eq:norm-PL}, we therefore obtain
\[
\dot{\Delta}(t)
\le
-\sqrt{2\mu_*\Delta(t)}
\qquad\text{for a.e. }t.
\]

On every interval where \(\Delta(t)>0\), we can divide by
\(2\sqrt{\Delta(t)}\) and obtain
\[
\frac{d}{dt}\sqrt{\Delta(t)}
=
\frac{\dot{\Delta}(t)}{2\sqrt{\Delta(t)}}
\le
-\sqrt{\frac{\mu_*}{2}}
\qquad\text{for a.e. }t.
\]
Integrating from \(0\) to \(t\) gives
\begin{equation}
\sqrt{\Delta(t)}
\le
\left[
\sqrt{\Delta(0)}
-
\sqrt{\frac{\mu_*}{2}}\,t
\right]_+,
\label{eq:finite-time-flow}
\end{equation}
where \([a]_+:=\max\{a,0\}\). Consequently, the trajectory reaches the
optimal level set in finite time, no later than
\begin{equation}
T_*
\le
\sqrt{\frac{2\Delta(0)}{\mu_*}}.
\label{eq:finite-time-T}
\end{equation}
Once \(\Delta(t)\) reaches zero, the monotonicity implied
by~\eqref{eq:flow-energy}, together with \(f(x(t))\ge f^\star\), gives
\[
f(x(t))=f^\star
\qquad\text{for all later times}.
\]

For the \(\ell_\infty\) sign flow, Ma et al.~\cite{ma2022dynamic} proved this finite-time phenomenon.
Equations~\eqref{eq:finite-time-flow}--\eqref{eq:finite-time-T} express the
same mechanism for an arbitrary norm and for every admissible velocity on a
nonunique optimal face.

Finally, strong convexity with respect to \(\|\cdot\|\) implies
the norm-PL inequality~\eqref{eq:norm-PL} with the same constant.
Indeed, if \(f\) is \(\mu_*\)-strongly convex with respect to
\(\|\cdot\|\), then for a minimizer \(x^\star\),
\[
f^\star
\ge
f(x)
+
\langle\nabla f(x),x^\star-x\rangle
+
\frac{\mu_*}{2}\|x^\star-x\|^2.
\]
Hence, by the dual-norm inequality,
\[
f(x)-f^\star
\le
\|\nabla f(x)\|_*\,\|x-x^\star\|
-
\frac{\mu_*}{2}\|x-x^\star\|^2
\le
\frac{\|\nabla f(x)\|_*^2}{2\mu_*},
\]
which is exactly~\eqref{eq:norm-PL}.

The canonical dynamics are already well defined as differential inclusions.
Filippov theory becomes useful when we choose a single-valued branch of these
inclusions, since such selectors can be discontinuous at zeros and ties.
We study this connection next.


\section{Filippov Theory and Sliding Updates}
\label{sec:filippov}

\paragraph{Why Filippov?}
The optimal-velocity maps in Section~\ref{sec:flow-geometry} are set-valued at zeros and ties. A numerical method, however, usually chooses a single direction. Any such selector is discontinuous when the chosen face changes. Filippov regularization gives a precise way to interpret the resulting discontinuous ODE and to recover the velocities that can appear at a switch.

\begin{remark}[Sliding and chattering]
Filippov regularization provides a rigorous way to describe trajectories at
discontinuities and, in particular, allows trajectories to slide along a
switching manifold. It does not, however, rule out all complicated switching
behavior. We say that a trajectory \emph{slides} when it remains on a switching
manifold for a nontrivial time interval. By \emph{chattering}, we mean the
occurrence of infinitely many switches in a finite time interval; the classical
Fuller phenomenon is a well-known example~\cite{ZelikinBorisov1994}.
We mention this distinction because Filippov regularization should not be viewed
as a general mechanism for eliminating chattering. In the analysis below, we
only use the existence and local sliding properties of Filippov solutions, and
no assumption excluding chattering is required.
\end{remark}

\subsection{Filippov Regularization in a Nutshell}
Let $F:\R^d\to\R^d$ be (possibly) discontinuous. Its \emph{Filippov set} at $x$ is
\begin{equation}
\mathcal{F}[F](x)
:= \bigcap_{\delta>0}\ \bigcap_{\substack{N\subset\R^d\\ m(N)=0}}
\overline{\operatorname{conv}}\{\,F(y):\ y\in B(x,\delta)\setminus N\,\},
\label{eq:filippov-closure}
\end{equation}
where $m$ denotes Lebesgue measure. The removal of arbitrary null sets is part of the classical Filippov definition~\cite{Filippov1988,Cortes2008Discontinuous}. Without it, one obtains a different local convexification in general.
A \emph{Filippov solution} is an absolutely continuous curve $x(\cdot)$ that satisfies
$\dot x(t)\in \mathcal{F}[F](x(t))$ for almost every $t$. Off switching sets (where $F$ is
continuous), $\mathcal{F}[F](x)=\{F(x)\}$ and Filippov reduces to the classical ODE.

\paragraph{Canonical sign flows and discontinuous selectors.}
The optimization dynamics themselves are the canonical inclusions from~\eqref{eq:canonical-flow}. In the two geometries used throughout the paper,
\begin{align}
\mathcal V_\infty(x)&=-\Sign(\nabla f(x)),\label{eq:Vinf}\\
\mathcal V_1(x)&=-\partial\|\cdot\|_\infty(\nabla f(x))
=\operatorname{conv}\{-s_i e^{(i)}:\ i\in\mathcal I(x),\ s_i\in\Sign(\partial_i f(x))\},
\label{eq:Vone}
\end{align}
with the usual interpretation $\mathcal V_1(x)=\{v:\|v\|_1\le1\}$ when $\nabla f(x)=0$.

A \emph{selector} of a set-valued map $\mathcal V$ is a single-valued map $b$ such that $b(x)\in\mathcal V(x)$ for every $x$. A continuous selector need not exist across a switching set. To connect the canonical inclusions with classical Filippov theory, choose measurable selectors
$b_\infty(x)\in\mathcal V_\infty(x)$ and $b_1(x)\in\mathcal V_1(x)$. Away from a switch the direction is unique. At a zero or tie, a fixed deterministic rule can choose one admissible direction. The following observation links the two descriptions.

\begin{proposition}[Filippov selectors remain in the optimal-velocity set]
\label{prop:selector-contained}
Let $b(x)\in\mathcal V(x)$ be any measurable selector of~\eqref{eq:canonical-V}. Then
\begin{equation}
\mathcal F[b](x)\subseteq\mathcal V(x)\qquad\text{for every }x.
\label{eq:selector-contained}
\end{equation}
Consequently, every Filippov solution generated by such a selector is also a solution of the canonical norm-constrained inclusion.
\end{proposition}
\begin{proof}
The map $\mathcal V$ is outer semicontinuous with compact convex values by Proposition~\ref{prop:canonical-existence}. Hence, for every $\varepsilon>0$, all values $\mathcal V(y)$ with $y$ sufficiently close to $x$ lie in $\mathcal V(x)+\varepsilon B$. The same is true for the selected vectors $b(y)$. Taking local closed convex hulls, removing null sets, and then letting $\varepsilon\downarrow0$ gives~\eqref{eq:selector-contained}.
\end{proof}

\begin{proposition}[Existence of solutions for the sign flows]\label{prop:filippov-existence}
Under Assumption~\ref{ass:standing}, the canonical sign maps $\mathcal V_\infty$ and $\mathcal V_1$ are Marchaud maps. Consequently, the inclusions $\dot x\in\mathcal V_\infty(x)$ and $\dot x\in\mathcal V_1(x)$ admit absolutely continuous solutions from every initial condition. Moreover, if $b_\infty$ and $b_1$ are measurable selectors of these maps, then the Filippov inclusions
$\dot x\in\mathcal F[b_\infty](x)$ and $\dot x\in\mathcal F[b_1](x)$ also admit absolutely continuous solutions from every initial condition; every such Filippov solution is a trajectory of the corresponding canonical inclusion.
\end{proposition}
\begin{proof}
The first statement is Proposition~\ref{prop:canonical-existence} specialized to the $\ell_\infty$ and $\ell_1$ norms. The selectors are measurable and globally bounded by the relevant unit ball. Standard Filippov existence results for measurable, locally essentially bounded vector fields therefore apply~\cite{Filippov1988,Cortes2008Discontinuous}. Proposition~\ref{prop:selector-contained} gives the final inclusion between the resulting trajectories.
\end{proof}

\subsection{Geometry of the Filippov Set on Tie Facets (\texorpdfstring{$\ell_1$}{l1} Flow)}

We will use the following local notion. An extreme velocity $v$ is \emph{essentially locally attainable} by a selector $b$ at $x$ if, for every $\delta>0$ and every $\varepsilon>0$, the set
\[
\{y\in B(x,\delta):\ \|b(y)-v\|<\varepsilon\}
\]
has positive Lebesgue measure. This is the natural condition because the Filippov construction ignores arbitrary null sets.

\begin{lemma}[Filippov set on tie facets for the $\ell_1$ flow]\label{lemma:FS}
Let $f\in C^1$ and
\[
\mathcal I(x):=\{i:\ |\partial_i f(x)|=\max_j|\partial_j f(x)|\}.
\]
Let $b_1$ be a measurable selector that, away from ties, chooses the unique direction
$-\operatorname{sign}(\partial_i f(x))e^{(i)}$. Then
\begin{equation}
\mathcal F[b_1](x)\subseteq\mathcal V_1(x)
=\operatorname{conv}\{-s_i e^{(i)}:\ i\in\mathcal I(x),\ s_i\in\Sign(\partial_i f(x))\}.
\label{eq:l1-filippov-inclusion}
\end{equation}
If every extreme point of $\mathcal V_1(x)$ is essentially locally attainable by the selector at $x$, then equality holds.
\end{lemma}
\begin{proof}
The inclusion is Proposition~\ref{prop:selector-contained}. For equality, essential local attainability puts every extreme direction of $\mathcal V_1(x)$ in every local essential closed convex hull entering the Filippov definition. Their convex hull is therefore contained in $\mathcal F[b_1](x)$, and the reverse inclusion follows from~\eqref{eq:l1-filippov-inclusion}.
\end{proof}

The attainability condition is not automatic at degenerate ties. For example, if $f(x_1,x_2)=\tfrac12x_2^2$, then at the origin both coordinates belong to $\mathcal I(0)$, but index~1 is never a maximizer in any punctured neighborhood. A selector based on the local maximizer therefore generates only the vertical segment, whereas the canonical set $\mathcal V_1(0)$ is the full $\ell_1$ unit ball. This is why we distinguish the optimization inclusion from the Filippov set of a particular selector.


\subsection{Equivalence off Switching Sets and Sliding on Them}

\begin{theorem}[a.e. equivalence and sliding]\label{thm:ae-equivalence}
Let $f\in C^1$. Define
\begin{align*}
\Sigma_\infty&:=\{x:\exists i,\ \partial_i f(x)=0\},\\
\Sigma_1&:=\{x:\nabla f(x)=0\ \text{or}\ \exists i\neq j:\
|\partial_i f(x)|=|\partial_j f(x)|=\|\nabla f(x)\|_\infty\}.
\end{align*}
Let $b_\infty$ and $b_1$ be measurable selectors of $\mathcal V_\infty$ and $\mathcal V_1$ that agree with the usual unique directions away from the switching sets. Then:
\begin{itemize}
\item[(i)] \textbf{Off switching sets.} If $x_0\notin\Sigma_\infty$, the signs of all partial derivatives are locally constant, so $b_\infty$ is locally constant and
$\mathcal F[b_\infty](x)=\{b_\infty(x)\}$ nearby. If $x_0\notin\Sigma_1$ and $\nabla f(x_0)\neq0$, the maximizing coordinate is locally unique with fixed sign, and the same conclusion holds for $b_1$. Thus the Filippov dynamics coincide a.e. with the ordinary selector dynamics whenever the trajectory remains off the switching set.
\item[(ii)] \textbf{On switching sets.} At every $x$,
\[
\mathcal F[b_\infty](x)\subseteq\mathcal V_\infty(x),\qquad
\mathcal F[b_1](x)\subseteq\mathcal V_1(x).
\]
Hence a Filippov trajectory on a persistent switching set uses convex combinations of locally attainable norm-steepest directions. Equality with the full canonical face holds whenever all of its extreme points are locally essentially attainable by the chosen selector. If the switching set is locally a $C^1$ manifold and the trajectory stays on it over an interval, its velocity is tangent to that manifold for a.e. time in the interval; this is the sliding case.
\end{itemize}
\end{theorem}
\begin{proof}
Part (i) follows from continuity of the partial derivatives and the positive gap between the unique active coordinate and the others. Part (ii) is Proposition~\ref{prop:selector-contained}, together with Lemma~\ref{lemma:FS} for the $\ell_1$ geometry. The equality statement follows from the same extreme-point argument used in Lemma~\ref{lemma:FS}. Finally, if a $C^1$ switching manifold is written locally as $\{h=0\}$ and $h(x(t))=0$ on an interval, the chain rule gives $\langle\nabla h(x(t)),\dot x(t)\rangle=0$ for a.e. $t$ on that interval.
\end{proof}

At degenerate switching points the inclusion can be strict. At a regular switching hypersurface, the Filippov set has an explicit form and gives a simple crossing-versus-sliding criterion.

\begin{proposition}[Regular switching surface: crossing or sliding]
\label{prop:regular-switch-slide}
Assume $f\in C^2$ near $\bar x$. Fix an index $i$ such that
$h(x):=\partial_i f(x)$ satisfies $h(\bar x)=0$ and $\nabla h(\bar x)\neq0$, while $\partial_j f(\bar x)\neq0$ for every $j\neq i$. By continuity, all signs except the $i$th one are fixed in a sufficiently small neighborhood of $\bar x$. Denote by $b^+$ and $b^-$ the two constant sign velocities on the open sides $h>0$ and $h<0$, respectively. Then
\begin{equation}
\mathcal F[b_\infty](\bar x)=\operatorname{conv}\{b^+,b^-\}.
\label{eq:regular-F}
\end{equation}
Moreover, a velocity tangent to the switching manifold $M=\{x:h(x)=0\}$ exists in this segment if and only if
\begin{equation}
\bigl\langle\nabla h(\bar x),b^+\bigr\rangle
\bigl\langle\nabla h(\bar x),b^-\bigr\rangle\le0.
\label{eq:slide-criterion}
\end{equation}
With strict inequality, the tangent Filippov velocity is unique. If
$\langle\nabla h,b^+\rangle<0<\langle\nabla h,b^-\rangle$, both side fields point toward $M$ and the sliding mode is attracting. The opposite strict ordering gives a repelling sliding mode. If the two normal components have the same nonzero sign, no convex combination is tangent and trajectories cross the surface locally.
\end{proposition}
\begin{proof}
By the implicit-function theorem, $M$ is a $C^1$ hypersurface near $\bar x$ and the two signs of $h$ occur on the two open sides of $M$. The selector is constant on each side, so the classical Filippov set is exactly the segment joining the two side limits, which gives~\eqref{eq:regular-F}. A vector $v=(1-\lambda)b^-+\lambda b^+$ is tangent to $M$ precisely when $\langle\nabla h(\bar x),v\rangle=0$. Such a $\lambda\in[0,1]$ exists exactly when the two endpoint normal components have opposite signs or one vanishes. Strict opposition gives a unique $\lambda$.
\end{proof}

\paragraph{Separable example.}
Suppose $f(x)=\sum_i f_i(x_i)$, each $f_i$ has a unique minimizer $x_i^\star$, and
$\operatorname{sign}(f_i'(x_i))=\operatorname{sign}(x_i-x_i^\star)$ whenever $x_i\neq x_i^\star$. Then the $\ell_\infty$ sign flow moves every nonoptimal coordinate toward its minimizer at unit speed. It reaches $x^\star$ after exactly $\|x(0)-x^\star\|_\infty$ time units. This simple example illustrates the finite-time behavior in~\eqref{eq:finite-time-flow}.

\subsection{Crossing and Sliding on a Simple Example}




The criterion above is easy to see in two dimensions. A trajectory may cross a switching manifold, or it may reach the manifold and then move along it.

These cases can be illustrated with
\[
f(x)=x_2+(x_2-a x_1)^2
\]
in a neighbourhood of the origin, where $a>0$. Let
\[
h(x):=x_2-a x_1,
\qquad
M=\{x:h(x)=0\}=\{x:x_2=a x_1\}.
\]
In a sufficiently small neighbourhood of $x=0$, we have
$1+2(x_2-a x_1)>0$ on both sides of $M$. Hence the velocity of the
$\ell_\infty$ sign flow $\dot x=-\sign(\nabla f)$ is
\[
\dot x=
\begin{cases}
b^+=(1,-1)^\top, & x_2>a x_1,\\
b^-=(-1,-1)^\top, & x_2<a x_1.
\end{cases}
\]

Since $\nabla h=(-a,1)^\top$,
\[
\langle \nabla h,b^+\rangle=-(a+1)<0,
\qquad
\langle \nabla h,b^-\rangle=a-1.
\]
If $a<1$, both quantities are negative, so the flow crosses the manifold
from the side $h>0$ to the side $h<0$. If $a>1$, the two quantities have
opposite signs: $b^+$ points toward $M$ from the side $h>0$, while $b^-$
points toward $M$ from the side $h<0$. The manifold is therefore attracting,
and the Filippov convexification contains a tangential sliding velocity.

Figure~\ref{fig:switching_sliding} illustrates these two cases: $a<1$ on the
left and $a>1$ on the right.

\begin{figure}[ht!]
\centering
\begin{tikzpicture}[x=2.15cm,y=2.15cm,>=Latex]
  \tikzset{
    axis/.style={->,line width=0.45pt},
    manifold/.style={line width=1.4pt},
    field/.style={->,line width=1.15pt},
    slide/.style={->,line width=1.8pt}
  }

  \begin{scope}[xshift=-2.25cm]
    \draw[axis] (-1.1,0)--(1.1,0) node[right] {$x_1$};
    \draw[axis] (0,-0.95)--(0,0.95) node[above] {$x_2$};

    \draw[manifold] (-1,-0.50)--(1,0.50)
      node[pos=0.97,above right] {$M$};

    \draw[field] (-0.82,0.66)--(-0.34,0.18)
      node[pos=0.55,above left,fill=white,inner sep=1pt] {$b^+$};

    \draw[field] (0.22,-0.06)--(-0.26,-0.54)
      node[pos=0.55,below right,fill=white,inner sep=1pt] {$b^-$};

    \node[align=center] at (0,-1.08) {crossing\\[-1pt]$a<1$};
  \end{scope}

  \begin{scope}[xshift=2.60cm]
    \draw[axis] (-1.1,0)--(1.1,0) node[right] {$x_1$};
    \draw[axis] (0,-0.95)--(0,0.95) node[above] {$x_2$};

    \draw[manifold] (-0.48,-0.82)--(0.48,0.82)
      node[pos=0.92,above right] {$M$};

    \draw[field] (-0.86,0.58)--(-0.42,0.14)
      node[pos=0.45,above left,fill=white,inner sep=1pt] {$b^+$};

    \draw[field] (0.80,0.10)--(0.36,-0.34)
      node[pos=0.45,below right,fill=white,inner sep=1pt] {$b^-$};

    \draw[slide] (0.16,0.28)--(-0.18,-0.30)
      node[pos=0.66,left,fill=white,inner sep=1pt] {$v_F$};

    \node[align=center] at (0,-1.08) {attracting sliding\\[-1pt]$a>1$};
  \end{scope}
\end{tikzpicture}
\caption{Crossing and sliding for $f(x)=x_2+(x_2-a x_1)^2$. The switching manifold is $M=\{x_2=a x_1\}$. When $a<1$, the two side velocities cross the manifold. When $a>1$, the side velocities point toward the manifold, and their Filippov convex hull contains the tangential sliding velocity $v_F$}
\label{fig:switching_sliding}
\end{figure}

\subsection{From Filippov to Practical Face-Aware Updates}

\paragraph{Two-hit sliding-track.}
We now translate this picture into a discrete rule. Near a crossing surface, the iteration should simply cross it. Near a sliding surface, repeated full sign steps may instead jump from one side to the other. We therefore use two consecutive \emph{strict} sign crossings of a partial derivative $\partial_i f$ as a local sliding indicator. For the three recorded values below, this means $d_{i,k-2}d_{i,k-1}<0$ and $d_{i,k-1}d_{i,k}<0$. When this happens, we approximate the condition $\partial_i f=0$ and damp the next coordinate velocity. 

More precisely, we take the standard step
\[
\vx_{k+1}=\vx_k-\eta_k\,\operatorname{sign}\!\big(\nabla f(\vx_k)\big),
\]
until two consecutive strict sign crossings of some partial derivative $\partial_i f$ are observed at steps $k-1$ and $k$. Let the corresponding values be $d_{i,k-2},d_{i,k-1},d_{i,k}$. They have alternating nonzero signs, and steps of magnitude $\eta_{k-2}$ and $\eta_{k-1}$ in coordinate $i$ have been performed between the three observations. We use a simple local model. We represent the contribution of all coordinates other than $i$ by $\alpha$ times the step size, and the contribution of coordinate $i$ by $\beta$ times the step size. This gives
\[ d_{i,k-1}-d_{i,k-2} = (\alpha + \beta)\eta_{k-2}, \qquad d_{i,k}-d_{i,k-1} = (\alpha - \beta)\eta_{k-1}.
\]
To target the manifold $\partial_i f=0$, we choose the next coordinate velocity so that
\[ d_{i,k+1}-d_{i,k} = -d_{i,k} = (\alpha + \beta\xi)\eta_k.
\]
Here $\xi$ controls the next velocity in coordinate $i$. When $0\le\xi\le1$, it damps the ordinary sign velocity while remaining in the admissible interval between the neighboring extreme sign velocities.

Solving the first two equations gives
\begin{equation}\label{eq:para_proj_ab}
\begin{aligned}
\alpha
&=\frac{d_{i,k-1}-d_{i,k-2}}{2\eta_{k-2}}
 +\frac{d_{i,k}-d_{i,k-1}}{2\eta_{k-1}},\\
\beta
&=\frac{d_{i,k-1}-d_{i,k-2}}{2\eta_{k-2}}
 -\frac{d_{i,k}-d_{i,k-1}}{2\eta_{k-1}}.
\end{aligned}
\end{equation}
For compactness, define
\begin{equation}\label{eq:D}
D:=d_{i,k}\eta_{k-2}
-d_{i,k-1}(\eta_{k-2}+\eta_{k-1})
+d_{i,k-2}\eta_{k-1}.
\end{equation}
Then solving the third equation for $\xi$ gives
\begin{equation}\label{eq:para_proj}
\xi=
\frac{
\begin{aligned}
&d_{i,k}\eta_k\eta_{k-2}
+d_{i,k-1}\eta_k\eta_{k-1}
+2d_{i,k}\eta_{k-1}\eta_{k-2}\\[-1mm]
&\qquad
-d_{i,k-1}\eta_k\eta_{k-2}
-d_{i,k-2}\eta_k\eta_{k-1}
\end{aligned}
}{\eta_k D}.
\end{equation}
The denominator in the expression for $\xi$ is therefore exactly $\eta_kD$.
If $D=0$ we regard the model as degenerate and keep the default sign step (switching regime).
In the equal-steps case $\eta_{k-2}=\eta_{k-1}=\eta_k$ this simplifies to
\[
D \;=\; \eta_k\,(d_{i,k}-2d_{i,k-1}+d_{i,k-2}),
\qquad
\xi \;=\; \frac{3d_{i,k}-d_{i,k-2}}{\,d_{i,k}-2d_{i,k-1}+d_{i,k-2}\,}.
\]

If $\xi>1$, the local model asks for a velocity beyond the admissible range and we take the ordinary sign step on that coordinate. If $\xi<0$, we clip it to zero. Thus the practical sliding correction uses only $\xi\in[0,1]$.

\paragraph{Convex-combination sliding on tie facets (for the $\ell_1$ flow).}
At $x_k$, let
\[
\mathcal I(x_k)=\Big\{i:\ |\partial_i f(x_k)|=\|\nabla f(x_k)\|_\infty\Big\}.
\]
The canonical $\ell_1$ velocity set is the convex hull of the corresponding signed basis directions. For a particular discontinuous selector, its Filippov set is contained in this hull and equals it under the local attainability condition of Lemma~\ref{lemma:FS}.

\medskip

\begin{lemma}[All convex combinations are maximally descending]\label{lemma:cvx_comb}
    Let $P=\max_j|\partial_j f(x)|$ and, for $i\in\mathcal I(x)$, pick $s_i\in\Sign(\partial_i f(x))$ and set
    $v_i=-s_i\,\ve^{(i)}$. Then for any convex weights $\{\alpha_i\}_{i\in\mathcal I}$
\[
v=\sum_{i\in\mathcal I}\alpha_i v_i
\quad\Rightarrow\quad
\langle \nabla f(x), v\rangle
= -\sum_{i\in\mathcal I}\alpha_i\,|\partial_i f(x)|\;=\;-P.
\]
\end{lemma}
\begin{proof}
  Linearity of the inner product and $|\partial_i f(x)|=P$ for $i\in\mathcal I$ (the case $P=0$ is trivial).
\end{proof}


\noindent\emph{The remainder of this subsection concerns the $\ell_1$ flow on tie facets.}

\paragraph{Discrete convex-combination update (tie-aware step).}
Given $\eta_k>0$ and active set $\mathcal I=\mathcal I(\vx_k)$, define signed coordinate moves
\[
\vx^{(i)}=\vx_k-\eta_k\,\operatorname{sign}\!\big(\partial_i f(\vx_k)\big)\,\ve^{(i)}\quad (i\in\mathcal I).
\]
Then update by any convex blend
\begin{equation}
\boxed{\quad
\vx_{k+1} \;=\; \sum_{i\in\mathcal I}\alpha_i\,\vx^{(i)},\qquad \alpha_i\ge 0,\ \sum_{i\in\mathcal I}\alpha_i=1.
\quad}
\label{eq:cc-update}
\end{equation}
Choices include: (a) \emph{vertex selection} (pick one $i^\star$ and set $\alpha_{i^\star}{=}1$), (b) \emph{equal weights} $\alpha_i=1/|\mathcal I|$, and (c) \emph{problem-driven weights} (for example, weights chosen to remain on a specific tie manifold). By Lemma~\ref{lemma:cvx_comb}, all these choices share the same first-order decrease
$\langle \nabla f(\vx_k), \vx_{k+1}-\vx_k\rangle=-\eta_k P$.

Figure~\ref{fig:l1-sliding-2d} illustrates the Filippov convexification for the $\ell_1$-constrained flow in two dimensions.

\begin{figure}[ht!]
\centering
\begin{tikzpicture}[scale=2, >=Latex]
  \draw[->] (-1.2,0) -- (1.2,0) node[below right] {$v_1$};
  \draw[->] (0,-1.2) -- (0,1.2) node[above left] {$v_2$};

  \draw[thick] (1,0) -- (0,1) -- (-1,0) -- (0,-1) -- cycle;

  \draw[->] (0,0) -- (0.8,0.8) node[above right] {$\nabla f$};

  \filldraw (-1,0) circle (0.03) node[below left] {$-\ve^{(1)}$};
  \filldraw (0,-1) circle (0.03) node[below right] {$-\ve^{(2)}$};

  \draw[line width=2pt, opacity=0.6] (-1,0) -- (0,-1);

  \draw[->, thick] (0,0) -- (-0.5,-0.5) node[below left] {$v^\star$};

  \node[align=center] at (-0.58,-0.12) {\scriptsize tie facet};
  \node[align=left] at (-0.6,-0.9) {\scriptsize $\mathcal V_1(x)$};
\end{tikzpicture}
\caption{Tie face for the $\ell_1$-constrained flow in 2D (velocity space).
The unit $\ell_1$ ball is a diamond. When $|\partial_1 f|=|\partial_2 f|$ (and, for illustration, both partial derivatives are positive), the active extreme directions are $-\ve^{(1)}$ and $-\ve^{(2)}$. Their convex hull (bold edge) is the canonical optimal-velocity face; a selector's Filippov set is contained in this edge and equals it when both directions are locally attainable. Any convex combination $v=\alpha(-\ve^{(1)})+(1-\alpha)(-\ve^{(2)})$ achieves the same instantaneous decrease $\langle \nabla f, v\rangle=-\|\nabla f\|_\infty$. The figure illustrates one possible sliding direction $v^\star$ (here, for $\alpha=1/2$); the specific sliding vector realized by a Filippov solution depends on higher-order properties of $f$}
\label{fig:l1-sliding-2d}
\end{figure}

\paragraph{Relation to the dual-norm formulation.}
For the $\ell_1$-constrained flow, the extreme directions on a tie facet are the vertices of the $\ell_1$ ball's face; their convex hull equals the subdifferential of $\|\cdot\|_\infty$ at $\nabla f$. Thus~\eqref{eq:cc-update} selects any element of $-\partial\|\cdot\|_\infty(\nabla f(\vx_k))$, consistent with the unifying Lemma~\ref{lemma:unified} in Section~\ref{subsec:dual-norm-lens}.

\subsection{Descent and Convergence for the $\ell_1$ Sliding Step}
Assume first that, for some $L>0$, $f$ satisfies the Euclidean quadratic upper bound
\[
f(y)\le f(x)+\langle\nabla f(x),y-x\rangle+\frac{L}{2}\|y-x\|_2^2.
\]
Let $P_k=\|\nabla f(\vx_k)\|_\infty=\max_i |\partial_i f(\vx_k)|$ and let
\[
\vg_k=\sum_{i\in\mathcal I(\vx_k)}\alpha_i\,\operatorname{sign}\!\big(\partial_i f(\vx_k)\big)\,\ve^{(i)}
\]
be the unit-$\ell_1$ search direction used in~\eqref{eq:cc-update}. By the upper bound above,
\[
f(\vx_{k+1}) \le f(\vx_k) + \langle \nabla f(\vx_k), -\eta_k \vg_k\rangle + \tfrac{L}{2}\,\eta_k^2 \|\vg_k\|_2^2.
\]

By Lemma~\ref{lemma:cvx_comb}, $\langle \nabla f(\vx_k), \vg_k\rangle=P_k$, hence
\begin{equation}
f(\vx_{k+1}) \;\le\; f(\vx_k)\;-\;\eta_k P_k\;+\;\frac{L}{2}\,\eta_k^2\,\|\vg_k\|_2^2.
\label{eq:cc-descent}
\end{equation}
Since $\vg_k$ has entries $\{\alpha_i\}_{i\in\mathcal I(\vx_k)}$ (up to signs), we have
$\|\vg_k\|_2^2=\sum_{i\in\mathcal I(\vx_k)}\alpha_i^2\le 1$.
Thus the uniform choice $\eta_k=P_k/L$ guarantees
\[
f(\vx_{k+1}) \le f(\vx_k) - \frac{P_k^2}{2L}.
\]
If the weights are known, the quadratic bound in~\eqref{eq:cc-descent} is minimized by the larger face-dependent step
$\eta_k=P_k/(L\|\vg_k\|_2^2)$. For equal weights on $r$ active coordinates, this becomes $\eta_k=rP_k/L$.

\paragraph{A convergence guarantee for tie-aware $\ell_1$ steps.}
The same argument can be written directly in the $\ell_1$ geometry. Assume that, for some $L_1>0$,
\begin{equation}
f(y)\le f(x)+\langle\nabla f(x),y-x\rangle+\frac{L_1}{2}\|y-x\|_1^2
\qquad\text{for all }x,y.
\label{eq:l1-smooth}
\end{equation}
Let $v_k$ be \emph{any} convex-combination velocity from~\eqref{eq:cc-update}, normalized so that $\|v_k\|_1=1$ when $\nabla f(x_k)\neq0$. Lemma~\ref{lemma:cvx_comb} gives
$\langle\nabla f(x_k),v_k\rangle=-\|\nabla f(x_k)\|_\infty$. Therefore the step
\begin{equation}
x_{k+1}=x_k+\frac{\|\nabla f(x_k)\|_\infty}{L_1}v_k
\label{eq:l1-scaled-slide}
\end{equation}
satisfies
\begin{equation}
f(x_{k+1})\le f(x_k)-\frac{\|\nabla f(x_k)\|_\infty^2}{2L_1}.
\label{eq:l1-slide-decrease}
\end{equation}
If, in addition,
\begin{equation}
\|\nabla f(x)\|_\infty^2\ge2\mu_1\bigl(f(x)-f^\star\bigr),
\label{eq:l1-pl}
\end{equation}
then every choice of convex weights obeys
\begin{equation}
f(x_{k+1})-f^\star\le\left(1-\frac{\mu_1}{L_1}\right)\bigl(f(x_k)-f^\star\bigr).
\label{eq:l1-slide-linear}
\end{equation}
Thus the freedom to slide on an active face does not alter the basic linear guarantee as long as the chosen velocity remains on the norm-steepest face. If $f$ is $\mu$-strongly convex in the Euclidean norm, then~\eqref{eq:l1-pl} holds with $\mu_1=\mu/d$ because $\|g\|_\infty\ge\|g\|_2/\sqrt d$.

\noindent\emph{Scope.} The convex-combination update above concerns the $\ell_1$ flow. The next two rules concern the $\ell_\infty$ sign flow and are independent of this tie-aware $\ell_1$ update.

\paragraph{Two face-aware SignGD rules.}
We consider two simple corrections near repeated sign changes:
\begin{itemize}
\item[(i)] \emph{One-hit freeze (freeze-on-first-flip).} If a coordinate crosses zero, we freeze that coordinate for the current step. This is the simplest correction.
\item[(ii)] \emph{Two-hit sliding-track.} To better approximate Filippov sliding, only after \emph{two consecutive strict sign crossings} on the same coordinate do we damp the $\pm1$ velocity on that coordinate, using the coefficient in~\eqref{eq:para_proj} to steer the partial derivative toward $0$.
\end{itemize}

\begin{algorithm}[H]
\caption{Face-aware SignGD (one-hit freeze)}
\label{alg:one_hit_freeze}
\begin{algorithmic}[1]
\STATE \textbf{Input:} $x_0$; positive stepsizes $\eta_k$ (for the guarantee below, $\eta_k=\|\nabla f(x_k)\|_1/L_\infty$)
\STATE Initialize $\texttt{prev\_g}\leftarrow \nabla f(x_0)$
\FOR{$k=0,1,2,\dots$}
  \STATE $g \leftarrow \nabla f(x_k)$,\quad choose $\eta_k$
  \STATE \textbf{Default step:} $x_{k+1} \leftarrow x_k - \eta_k\,\sign(g)$ \hfill (componentwise, $\sign(0)=0$)
  \STATE \textbf{Freeze on first crossing:} for any $i$ with $\texttt{prev\_g}_i\,g_i<0$, set $x_{k+1,i}\leftarrow x_{k,i}$
  \STATE $\texttt{prev\_g}\leftarrow g$
\ENDFOR
\end{algorithmic}
\end{algorithm}

\begin{algorithm}[H]
\caption{Face-aware SignGD (two-hit sliding-track)}
\label{alg:two_hit_sliding}
\begin{algorithmic}[1]
\STATE \textbf{Input:} $x_0$; positive stepsizes $\eta_k$ (same convention as above)
\STATE Initialize $\texttt{g\_pprev}\!\leftarrow \nabla f(x_0)$,\ $\texttt{g\_prev}\!\leftarrow \nabla f(x_0)$; keep the two most recent stepsizes once available
\FOR{$k=0,1,2,\dots$}
  \STATE $g \leftarrow \nabla f(x_k)$,\quad choose $\eta_k$
  \STATE \textbf{Default velocity:} $u \leftarrow -\sign(g)$ \hfill ($\|u\|_\infty\le 1$; $\sign(0)=0$)
  \IF{$k\ge 2$}
    \FOR{$i=1,\dots,d$}
      \STATE \textbf{Two-hit test:} if $g_i\,\texttt{g\_prev}_i<0$ and $\texttt{g\_prev}_i\,\texttt{g\_pprev}_i<0$ then
      \STATE \quad Let $d_{i,k-2}\!=\texttt{g\_pprev}_i$, $d_{i,k-1}\!=\texttt{g\_prev}_i$, $d_{i,k}\!=g_i$
      \STATE \quad Compute $D$ from~\eqref{eq:D}
      \STATE \quad If $D\neq 0$, set $\displaystyle \xi$ as in~\eqref{eq:para_proj}.
      \STATE \quad \textbf{Clamp and set:} if $D\neq 0$, set $u_i \leftarrow -\,\sign(g_i)\cdot \mathrm{clip}(\xi,0,1)$
    \ENDFOR
  \ENDIF
  \STATE $x_{k+1} \leftarrow x_k + \eta_k\,u$ \hfill (i.e., $x_{k+1}=x_k-\eta_k\sign(g)$ off sliding)
  \STATE Store $g$ and $\eta_k$ in the two-step history for the next iterations
\ENDFOR
\end{algorithmic}
\end{algorithm}

\begin{remark}[Equal-step simplification and safety]
If $\eta_{k-2}=\eta_{k-1}=\eta_k$, the formula simplifies to $\ \xi=\frac{3d_{i,k}-d_{i,k-2}}{d_{i,k}-2d_{i,k-1}+d_{i,k-2}}$.
Clamping $\xi$ to $[0,1]$ keeps $u_i$ in the admissible interval $[-1,1]$ and never reverses the current descent sign. If the denominator vanishes or $\xi>1$, the rule reduces to the default sign step on that coordinate.
\end{remark}

\subsection{A Convergence Safeguard for Face-Aware SignGD}
\label{subsec:safeguarded-sliding}

Both face-aware rules produce damped versions of the ordinary sign direction, which we now write in a common form. The one-hit rule avoids an immediate recrossing, while the two-hit rule tries to approximate a local sliding velocity. We then add a simple safeguard that checks whether the proposed direction retains enough of the descent guaranteed for ordinary SignGD.

Assume that $f$ is $L_\infty$-smooth with respect to the $\ell_\infty$ norm:
\begin{equation}
f(y)\le f(x)+\langle\nabla f(x),y-x\rangle+\frac{L_\infty}{2}\|y-x\|_\infty^2
\qquad\text{for all }x,y.
\label{eq:linf-smooth}
\end{equation}
Let $g=\nabla f(x)\neq0$. Any one-hit or two-hit candidate can be written as
\begin{equation}
d(w,g)=w\odot\operatorname{sign}(g),\qquad 0\le w_i\le1.
\label{eq:damped-sign-direction}
\end{equation}
For the one-hit rule, $w_i\in\{0,1\}$; for the two-hit rule, a detected coordinate uses the clipped coefficient $\xi_i$ from~\eqref{eq:para_proj}. Define
\begin{equation}
A(w,g):=\sum_i w_i|g_i|,\qquad M(w):=\|w\|_\infty,
\label{eq:A-M}
\end{equation}
and the dimensionless descent score
\begin{equation}
q(w,g):=\frac{2A(w,g)}{\|g\|_1}-M(w)^2.
\label{eq:quality-q}
\end{equation}
Since $A(w,g)\le M(w)\|g\|_1$ and $0\le M(w)\le1$, one has $-1\le q(w,g)\le1$. The ordinary sign direction corresponds to $w=\mathbf 1$ and has $q=1$. Notice that $\|d(w,g)\|_\infty\le M(w)$; equality can fail only when a largest weight is attached to a zero gradient coordinate. Using $M(w)$ therefore gives a conservative certificate and matches the quantity computed in the experiments.

\begin{lemma}[Best certified step along a damped sign direction]
\label{lem:damped-sign-optimal}
Under~\eqref{eq:linf-smooth}, let $d=d(w,g)$ and suppose $M(w)>0$. Then, for every $\eta\ge0$,
\[
f(x-\eta d)\le f(x)-\eta A(w,g)+\frac{L_\infty}{2}\eta^2M(w)^2.
\]
The right-hand side is minimized at
\[
\eta^\star=\frac{A(w,g)}{L_\infty M(w)^2},
\]
and therefore
\begin{equation}
f(x-\eta^\star d)\le f(x)-\frac{A(w,g)^2}{2L_\infty M(w)^2}.
\label{eq:damped-optimal-decrease}
\end{equation}
Equivalently, if
$r(w,g):=A(w,g)/(M(w)\|g\|_1)$, the guaranteed decrease is
$r(w,g)^2\|g\|_1^2/(2L_\infty)$. If $M(w)=0$, the candidate is the zero direction.
\end{lemma}
\begin{proof}
Apply~\eqref{eq:linf-smooth} to $y=x-\eta d$ and use
$\langle g,d\rangle=A(w,g)$ and $\|d\|_\infty\le M(w)$. The remaining upper bound is a scalar quadratic in $\eta$.
\end{proof}

The preceding lemma gives the best step certified by the $\ell_\infty$ model for any damped sign direction. For the two-hit rule, however, we prefer to keep the same nominal step as ordinary scaled SignGD, because that is the step used by the interpolation model. This leads to the following same-step certificate.

\begin{lemma}[Descent certificate for a face-aware sign candidate]
\label{lem:damped-sign}
Under~\eqref{eq:linf-smooth}, use the same scaled SignGD step length
\begin{equation}
\eta=\frac{\|g\|_1}{L_\infty},
\qquad x^+=x-\eta\,d(w,g).
\label{eq:damped-step}
\end{equation}
Then
\begin{equation}
f(x^+)\le f(x)-\frac{q(w,g)}{2L_\infty}\|g\|_1^2.
\label{eq:damped-decrease}
\end{equation}
\end{lemma}
\begin{proof}
Apply~\eqref{eq:linf-smooth} with $y=x-\eta d(w,g)$. We have
\[
\langle g,d(w,g)\rangle=A(w,g),
\qquad
\|d(w,g)\|_\infty\le M(w).
\]
Therefore,
\[
f(x-\eta d)\le f(x)-\eta A(w,g)+\frac{L_\infty}{2}\eta^2M(w)^2.
\]
Substituting~\eqref{eq:damped-step} gives
\[
f(x^+)\le f(x)-\frac{\|g\|_1^2}{2L_\infty}
\left(\frac{2A(w,g)}{\|g\|_1}-M(w)^2\right),
\]
which is~\eqref{eq:damped-decrease}.
\end{proof}

Fix a threshold $\rho\in(0,1]$. At iteration $k$, first generate the candidate weights $w_k$ by the one-hit or two-hit rule, using the same nominal step length $\eta_k=\|g_k\|_1/L_\infty$ that appears in the update and, for the two-hit rule, in the interpolation formula. Compute $q_k=q(w_k,g_k)$. If $q_k<\rho$, discard only the face correction and take the ordinary sign direction $w_k=\mathbf1$; otherwise keep the candidate. The detector, the two-hit construction, and its coefficient $\xi$ remain part of the algorithm; the safeguard only decides whether the proposed correction is sufficiently descent-preserving.

\begin{theorem}[Global rate of safeguarded face-aware SignGD]
\label{thm:safeguarded-sliding}
Assume~\eqref{eq:linf-smooth} and the dual-norm PL inequality
\begin{equation}
\|\nabla f(x)\|_1^2\ge2\mu_\infty\bigl(f(x)-f^\star\bigr)
\qquad\text{for all }x.
\label{eq:linf-pl}
\end{equation}
At every iteration, generate the one-hit or two-hit candidate, apply the safeguard $q_k\ge\rho$ described above, and use
$\eta_k=\|\nabla f(x_k)\|_1/L_\infty$. Then
\begin{equation}
f(x_{k+1})-f^\star
\le\left(1-\rho\frac{\mu_\infty}{L_\infty}\right)
\bigl(f(x_k)-f^\star\bigr),
\label{eq:safeguarded-rate}
\end{equation}
and consequently
\begin{equation}
f(x_k)-f^\star
\le\left(1-\rho\frac{\mu_\infty}{L_\infty}\right)^k
\bigl(f(x_0)-f^\star\bigr).
\label{eq:safeguarded-global}
\end{equation}
Euclidean $\mu$-strong convexity implies~\eqref{eq:linf-pl} with $\mu_\infty=\mu$.
\end{theorem}
\begin{proof}
After the safeguard, either the candidate is accepted with $q_k\ge\rho$ or it is replaced by the ordinary sign direction, for which $q_k=1$. Lemma~\ref{lem:damped-sign} therefore gives
\[
f(x_{k+1})\le f(x_k)-\frac{\rho}{2L_\infty}\|\nabla f(x_k)\|_1^2.
\]
Using~\eqref{eq:linf-pl} yields~\eqref{eq:safeguarded-rate}, and iteration gives~\eqref{eq:safeguarded-global}. Under Euclidean $\mu$-strong convexity,
$\|\nabla f(x)\|_1\ge\|\nabla f(x)\|_2\ge\sqrt{2\mu(f(x)-f^\star)}$.
\end{proof}

\begin{remark}[Connection with the separable-smoothness step]
Assumption~\ref{ass:coord} implies~\eqref{eq:linf-smooth} with
$L_\infty=\sum_iL_i=\|\bar L\|_1$, because
\[
\sum_iL_i(y_i-x_i)^2
\le
\|\bar L\|_1\|y-x\|_\infty^2.
\]
Thus the safeguarded method uses the same adaptive step as Section~\ref{sec:signgd-strong},
\[
\eta_k=\frac{\|\nabla f(x_k)\|_1}{\|\bar L\|_1},
\]
and applies the test~\eqref{eq:quality-q} to the face-aware candidate.
\end{remark}

\paragraph{Interpretation.}
The theorem does not assume that the two-hit model always identifies a true sliding manifold. When the local interpolation removes too much descent, the observable score $q_k$ rejects the correction and the method returns to the ordinary scaled sign step. Thus the geometric rule proposes the correction, and the safeguard decides whether to keep it.

In summary, Filippov theory explains how a discontinuous selector behaves at a switch. Off the switching set it agrees with the ordinary ODE. On a switching manifold it allows convexified velocities, including tangential sliding velocities when they are locally attainable. This directly motivates the one-hit, two-hit, and tie-aware updates above.

We now turn to the basic discrete SignGD method and its convergence under strong convexity.


\section{Deterministic Convergence of SignGD under Strong Convexity}
\label{sec:signgd-strong}

We study the basic SignGD iteration
\begin{equation}
  \vx_{k+1} = \vx_k - \eta_k\, \sign\!\big(\nabla f(\vx_k)\big),
  \qquad \eta_k>0,
  \label{eq:signgd-iter}
\end{equation}
under \assref{ass:coord} (separable smoothness) and, when stated, \assref{ass:strong} (strong convexity).

Linear convergence of scaled SignGD under related assumptions was established by Li et al.~\cite{li2021faster}, while the $\ell_\infty$ steepest-descent interpretation is discussed by Balles, Pedregosa, and Le Roux (arXiv 2002.08056, 2020). The short derivation below provides the estimates used in the face-aware refinements. We use the term \emph{separable smoothness} for Assumption~\ref{ass:coord}; diagonal Hessian entries alone do not certify this assumption for a generic dense objective.

Throughout this section we write

\begin{equation*}
\begin{aligned}
f^\star&:=\min_{x\in\mathbb{R}^d}f(x),
&\Delta_k&:=f(\vx_k)-f^\star,\\
s_k&:=\sign\!\bigl(\nabla f(\vx_k)\bigr)\in\{-1,0,+1\}^d,
&\|\bar L\|_1&:=\sum_{i=1}^d L_i,\qquad L:=\max_iL_i.
\end{aligned}
\end{equation*}
and we also use the active-face curvature
\[
S_k \;:=\; \sum_{i:\,\partial_i f(\vx_k)\neq 0} L_i \;\le\; \|\bar L\|_1
\] 
when deriving refinement bounds.

By \assref{ass:coord}, for any $\vx,\vy\in\R^d$ we have the separable quadratic upper bound
\begin{equation}
f(\vy) \;\le\; f(\vx) + \langle \nabla f(\vx), \vy-\vx\rangle
            + \tfrac12 \sum_{i=1}^d L_i\,(\vy_i-\vx_i)^2,
\label{eq:coord-smooth}
\end{equation}
which in particular implies the Euclidean quadratic upper bound with constant $L=\max_iL_i$, since
$\sum_i L_i(\vy_i-\vx_i)^2 \le L\,\|\vy-\vx\|_2^2$.

\subsection{Two Basic Inequalities}
We begin with a norm relation that connects $\Delta_k$ to gradient norms.

\begin{lemma}[Gradient--suboptimality bound]
\label{lem:norm-rel}
For any $k$,
\[
\|\nabla f(\vx_k)\|_1 \;\ge\; \|\nabla f(\vx_k)\|_2 \;\ge\; \sqrt{2\mu\,\Delta_k}.
\]
\end{lemma}

\begin{proof}
The inequality $\|\cdot\|_1\ge \|\cdot\|_2$ is standard. For the second, $\mu$--strong convexity
(\assref{ass:strong}) and optimality of $\vx^\star$ give
\[
f^\star \;\ge\; f(\vx_k) + \langle \nabla f(\vx_k), \vx^\star-\vx_k\rangle + \tfrac{\mu}{2}\|\vx^\star-\vx_k\|_2^2,
\]
hence
\(
\Delta_k \le \max_{r\ge0}\{\|\nabla f(\vx_k)\|_2\,r - \tfrac{\mu}{2}r^2\}
= \|\nabla f(\vx_k)\|_2^2/(2\mu).
\)
\end{proof}

\begin{lemma}[Per-iteration descent under \assref{ass:coord}]
\label{lem:descent}
One step of~\eqref{eq:signgd-iter} satisfies
\begin{equation}
\begin{aligned}
\Delta_{k+1}
&\le
\Delta_k-\eta_k\|\nabla f(\vx_k)\|_1
+\frac{\eta_k^2}{2}\sum_{i=1}^dL_i(s_{k,i})^2\\
&\le
\Delta_k-\eta_k\|\nabla f(\vx_k)\|_1
+\frac{\eta_k^2}{2}\|\bar L\|_1.
\end{aligned}
\label{eq:basic-descent}
\end{equation}
\end{lemma}

\begin{proof}
Apply the separable upper bound~\eqref{eq:coord-smooth} with $\vy=\vx_k-\eta_k s_k$ and $\vx=\vx_k$:
\[
f(\vx_{k+1})
\le f(\vx_k) + \langle \nabla f(\vx_k), -\eta_k s_k\rangle
 + \tfrac12\sum_{i=1}^d L_i (\eta_k s_{k,i})^2.
\]
Since $\langle \nabla f(\vx_k), s_k\rangle=\sum_i |\partial_i f(\vx_k)|=\|\nabla f(\vx_k)\|_1$
and $(s_{k,i})^2\le 1$, subtract $f^\star$ to obtain~\eqref{eq:basic-descent}.
\end{proof}

\subsection{A Sufficient Decrease Inequality and Its Minimizer}
Inequality~\ref{eq:basic-descent} is a one-step quadratic upper bound in $\eta_k$.
When $\nabla f(x_k)\neq0$, minimizing its right-hand side over $\eta\ge0$ yields the best step for that quadratic bound:
\[
\eta_k^{\rm quad}
\;=\;
\frac{\|\nabla f(\vx_k)\|_1}{\sum_{i=1}^d L_i (s_{k,i})^2}
\;=\;
\frac{\|\nabla f(\vx_k)\|_1}{S_k}
\;\;\ge\;\;
\frac{\|\nabla f(\vx_k)\|_1}{\|\bar L\|_1},
\]
where we used $S_k=\sum_{i:\,\partial_i f(\vx_k)\neq 0} L_i=\sum_i L_i(s_{k,i})^2\le \|\bar L\|_1$.

Using the conservative denominator $\|\bar L\|_1$ gives a simple, implementable rule and a direct decrease bound.

\begin{proposition}[Sufficient decrease with an adaptive step]
\label{prop:suff-decrease}
Under \assref{ass:coord}, with
\(
\displaystyle
\eta_k := \frac{\|\nabla f(\vx_k)\|_1}{\|\bar L\|_1},
\)
the SignGD step satisfies
\begin{equation}
\Delta_{k+1}
\;\le\;
\Delta_k \;-\; \frac{1}{2\|\bar L\|_1}\,\|\nabla f(\vx_k)\|_1^2.
\label{eq:suff-decrease}
\end{equation}
If, in addition, \assref{ass:strong} holds, then
\begin{equation}
\Delta_{k+1}
\;\le\;
\Bigl(1-\frac{\mu}{\|\bar L\|_1}\Bigr)\,\Delta_k.
\end{equation}
\end{proposition}

\begin{proof}
Plug $\eta_k$ into~\eqref{eq:basic-descent} to get
\(
\Delta_{k+1}
\le \Delta_k - \frac{\|\nabla f(\vx_k)\|_1^2}{\|\bar L\|_1}
     + \frac{1}{2}\frac{\|\nabla f(\vx_k)\|_1^2}{\|\bar L\|_1}
= \Delta_k - \frac{1}{2\|\bar L\|_1}\|\nabla f(\vx_k)\|_1^2,
\)
which is \eqref{eq:suff-decrease}. Under \assref{ass:strong}, Lemma~\ref{lem:norm-rel} gives
$\|\nabla f(\vx_k)\|_1^2 \ge 2\mu \Delta_k$, yielding the linear contraction.
\end{proof}

\paragraph{Discussion.}
Inequality~\eqref{eq:suff-decrease} guarantees monotone decrease of $f(\vx_k)$ and a per-step contraction governed by $\mu/\|\bar L\|_1$ under strong convexity. The denominator $\|\bar L\|_1=\sum_i L_i$ reflects the non-Euclidean (separable) quadratic model. When available, replacing $\|\bar L\|_1$ by the active-face curvature $S_k=\sum_{i:\,\partial_i f(\vx_k)\neq 0} L_i$ sharpens both the step and the decrease bound.

\subsection{Linear Convergence}
We now state the rate in function value and derive a corresponding distance decay.

\begin{theorem}[Linear rate with adaptive step]\label{thm:linear}
Under \assref{ass:coord} and \assref{ass:strong}, with the adaptive step
\(
\displaystyle
\eta_k = \frac{\|\nabla f(\vx_k)\|_1}{\|\bar L\|_1},
\)
the SignGD iterates~\eqref{eq:signgd-iter} obey
\begin{equation}
\Delta_k \;\le\; \Bigl(1-\frac{\mu}{\|\bar L\|_1}\Bigr)^k \Delta_0
\;\le\;
\Delta_0\,\exp\!\Bigl(-\frac{\mu}{\|\bar L\|_1}\,k\Bigr),
\label{eq:linear-rate}
\end{equation}
and, furthermore,
\begin{equation}
\|\vx_k-\vx^\star\|_2^2
\;\le\; \frac{L}{\mu}\,\Bigl(1-\frac{\mu}{\|\bar L\|_1}\Bigr)^k\,\|\vx_0-\vx^\star\|_2^2,
\qquad L=\max_i L_i.
\label{eq:distance-rate}
\end{equation}
Hence it is sufficient to take
\[
k \;\ge\; \frac{\|\bar L\|_1}{\mu}\,\log\!\Bigl(\frac{\Delta_0}{\varepsilon}\Bigr)
\]
to guarantee $f(\vx_k)-f^\star\le\varepsilon$.
\end{theorem}

\begin{proof}
By Proposition~\ref{prop:suff-decrease}, 
\(
\Delta_{k+1}\le \Delta_k - \tfrac{1}{2\|\bar L\|_1}\|\nabla f(\vx_k)\|_1^2.
\)
Under \assref{ass:strong}, Lemma~\ref{lem:norm-rel} yields 
\(\|\nabla f(\vx_k)\|_1^2 \ge 2\mu\,\Delta_k\), hence
\(\Delta_{k+1}\le (1-\mu/\|\bar L\|_1)\Delta_k\), which telescopes to \eqref{eq:linear-rate}.
For \eqref{eq:distance-rate}, combine strong convexity
\(\Delta_k \ge (\mu/2)\|\vx_k-\vx^\star\|_2^2\)
with the Euclidean quadratic upper bound from \assref{ass:coord} at $k=0$,
\(\Delta_0 \le (L/2)\|\vx_0-\vx^\star\|_2^2\), to obtain
\[
\|\vx_k-\vx^\star\|_2^2
\le \frac{2}{\mu}\,\Delta_k
\le \frac{2}{\mu}\,\Bigl(1-\frac{\mu}{\|\bar L\|_1}\Bigr)^k \Delta_0
\le \frac{L}{\mu}\,\Bigl(1-\frac{\mu}{\|\bar L\|_1}\Bigr)^k \|\vx_0-\vx^\star\|_2^2.
\]
\end{proof}

\noindent
If $\nabla f(\vx_k)=0$, then $\eta_k=0$ and $\vx_{k+1}=\vx_k$, so the recursion terminates at an optimizer.

\paragraph{Active-face refinement.}
At iteration $k$, only coordinates with non-zero gradient can move:
let $\mathcal I_k:=\{i:\ \partial_i f(\vx_k)\neq 0\}$ and
$S_k:=\sum_{i\in\mathcal I_k} L_i\le \|\bar L\|_1$.
If we replace the conservative step by the "face-aware" step
\(
\eta_k = \|\nabla f(\vx_k)\|_1 / S_k,
\)
then the proof of Proposition~\ref{prop:suff-decrease} yields
\(
\Delta_{k+1}\le \bigl(1-\mu/S_k\bigr)\Delta_k.
\)
When there exists $S_{\max}$ with $S_k\le S_{\max}$ for all $k$ (e.g., if $|\mathcal I_k|$ stays bounded and $L_i$ are comparable), we obtain the improved global rate
\(
\Delta_k\le (1-\mu/S_{\max})^k\Delta_0.
\)

The exact active set above is mainly useful when the model has coordinates whose derivatives become exactly zero. A thresholded or deliberately sparse face can also be analyzed, but then one must account for the gradient mass that is discarded.

\begin{proposition}[Mass-aware restriction]
\label{prop:mass-aware-face}
Assume separable smoothness and strong convexity. At $x_k$, let $J_k\subseteq\{1,\ldots,d\}$ be any nonempty set and define
\[
A_k:=\sum_{i\in J_k}|\partial_i f(x_k)|,
\qquad
S(J_k):=\sum_{i\in J_k}L_i,
\qquad
\theta_k:=\frac{A_k}{\|\nabla f(x_k)\|_1}\in[0,1].
\]
Let $d_{k,i}=\operatorname{sign}(\partial_i f(x_k))$ for $i\in J_k$ and $d_{k,i}=0$ otherwise, and take
\begin{equation}
x_{k+1}=x_k-\frac{A_k}{S(J_k)}d_k.
\label{eq:masked-face-step}
\end{equation}
Then
\begin{equation}
f(x_{k+1})\le f(x_k)-\frac{A_k^2}{2S(J_k)},
\label{eq:masked-face-descent}
\end{equation}
and
\begin{equation}
\Delta_{k+1}\le
\left(1-\theta_k^2\frac{\mu}{S(J_k)}\right)\Delta_k.
\label{eq:masked-face-rate}
\end{equation}
In particular, if $A_k\ge\theta\|\nabla f(x_k)\|_1$ and $S(J_k)\le \bar S$ uniformly, then
$\Delta_k\le(1-\theta^2\mu/\bar S)^k\Delta_0$.
\end{proposition}
\begin{proof}
Apply~\eqref{eq:coord-smooth} to~\eqref{eq:masked-face-step}. The linear term is
$-\eta A_k$ and the quadratic term is $\eta^2S(J_k)/2$. Choosing $\eta=A_k/S(J_k)$ gives~\eqref{eq:masked-face-descent}. Since
$A_k=\theta_k\|\nabla f(x_k)\|_1$ and Lemma~\ref{lem:norm-rel} gives
$\|\nabla f(x_k)\|_1^2\ge2\mu\Delta_k$, equation~\eqref{eq:masked-face-rate} follows.
\end{proof}

This proposition justifies the use of approximate active faces. A numerical threshold by itself is not enough: what matters for the rate is the retained gradient mass $\theta_k$, together with the curvature $S(J_k)$. The exact active-face refinement is recovered by taking $J_k=\mathcal I_k$, for which $\theta_k=1$.

\begin{proposition}[When $S_k \ll \|\bar L\|_1$]\label{prop:face-vs-Lbar}
Under Assumptions~\ref{ass:coord}--\ref{ass:strong}, let $r_k:=|\mathcal I_k|$, $L_{\min}:=\min_iL_i$, $L_{\max}:=\max_iL_i$, and $\kappa_L:=L_{\max}/L_{\min}$. Strong convexity implies $L_i\ge\mu>0$, so $\kappa_L$ is well defined. Then
\[
\frac{r_k}{d\,\kappa_L}\ \le\ \frac{S_k}{\|\bar L\|_1}\ \le\ \frac{r_k\,\kappa_L}{d}.
\]
In particular, if the $L_i$ are comparable (moderate $\kappa_L$) and only $r_k\!\ll d$ coordinates are active, then $S_k\ll \|\bar L\|_1$, so the face-aware contraction $1-\mu/S_k$ is strictly sharper than $1-\mu/\|\bar L\|_1$.
\end{proposition}

\begin{proof}
First, $L_i\ge\mu$ for every $i$: apply the separable upper model and the strong-convexity lower model to $y=x+t e^{(i)}$ and compare their quadratic terms. Hence $L_{\min}>0$. Now $r_k L_{\min}\le S_k\le r_k L_{\max}$ and $dL_{\min}\le\|\bar L\|_1\le dL_{\max}$; dividing the bounds gives the result.
\end{proof}

\begin{corollary}[Equal-curvature case]\label{cor:equal-L}
If $L_i\equiv L_0$, then $S_k/\|\bar L\|_1=r_k/d$. Thus, at iteration $k$, the active-coordinate denominator is smaller by the factor $r_k/d$. If $r_k\le r$ uniformly, the corresponding global bound uses $rL_0$ instead of $dL_0$, giving a factor $d/r$ improvement in the bound.
\end{corollary}

By Proposition~\ref{prop:face-vs-Lbar}, improvements are most pronounced when $r_k/d$ is small and the $L_i$ are of comparable size.

\paragraph{Comparison to Euclidean GD.}
Euclidean gradient descent with step $1/L$, using the same quadratic upper-model constant, contracts by $(1-\mu/L)$.
Our factor $(1-\mu/\|\bar L\|_1)$ is worse whenever $\|\bar L\|_1\gg L$,
which is typical if many coordinates contribute simultaneously.
The difference comes from the curvature model used by the sign step: the bound in~\eqref{eq:basic-descent} adds coordinate curvatures, whereas Euclidean GD uses a spectral smoothness constant. Sign updates are nevertheless attractive when their simple directions, communication cost, or robustness are important.

\paragraph{Implementation notes.}
The step rule $\eta_k=\|\nabla f(\vx_k)\|_1/\|\bar L\|_1$ is scale-free and requires only
(i) the gradient and (ii) a bound on $\|\bar L\|_1$ (e.g., known from model structure or estimated at initialization).
If a useful bound on $\|\bar L\|_1$ is unavailable, a standard backtracking test along the sign direction can be used instead. It preserves a sufficient-decrease condition without requiring an explicit global curvature constant.

While the basic scheme converges linearly, inertia can improve practical speed. 
We therefore study a safeguarded inertial variant and provide a rigorous descent guarantee.

\section{Inertial Sign Gradient Descent with Restart}
\label{sec:inertial-signgd}

The basic SignGD bound can be slow when $\|\bar L\|_1/\mu$ is large. This motivates adding an inertial extrapolation before the sign step. 

Momentum is also useful for normalized and sign methods in stochastic settings; see, for example,~\cite{pmlr-v119-cutkosky20b}. Here we only study the deterministic strongly convex case.

We therefore study a deterministic inertial variant of SignGD for $\mu$-strongly convex and separably smooth functions. We use a restart safeguard to preserve monotone descent when extrapolation is not useful. We do not claim an accelerated complexity bound: establishing a sharper rate for the inertial scheme remains an open question.


We consider the following inertial SignGD algorithm:

\begin{algorithm}[H]
\caption{Inertial Sign Gradient Descent with restart}
\label{alg:momentum-signgd}
\begin{algorithmic}[1]
\STATE \textbf{Input:} $x_0=x_1\in\mathbb R^d$, momentum parameter $\beta\in[0,1)$, step-size policy
\FOR{$k=1,2,\dots$}
    \STATE $v_k \leftarrow x_k+\beta(x_k-x_{k-1})$
    \STATE \textbf{Restart:} if $f(v_k)>f(x_k)$, set $v_k\leftarrow x_k$
    \STATE $g_k\leftarrow\nabla f(v_k)$ and choose $\eta_k>0$
    \STATE $x_{k+1}\leftarrow v_k-\eta_k\operatorname{sign}(g_k)$
\ENDFOR
\end{algorithmic}
\end{algorithm}

The restart is simple: extrapolation is kept only when it does not increase the objective. The experiments also report the same inertial update without restart. For the theoretical statement below we use the adaptive step
\begin{equation}
\eta_k=\frac{\|\nabla f(v_k)\|_1}{\|\bar L\|_1}.
\label{eq:inertial-adaptive-step}
\end{equation}

\begin{theorem}[Linear convergence of the restarted inertial scheme]
\label{thm:inertial-linear}
Assume separable smoothness (Assumption~\ref{ass:coord}) and Euclidean $\mu$-strong convexity (Assumption~\ref{ass:strong}). Let $(x_k)$ be generated by Algorithm~\ref{alg:momentum-signgd} with the adaptive step~\eqref{eq:inertial-adaptive-step}. Then
\begin{equation}
f(x_{k+1})-f^\star
\le
\left(1-\frac{\mu}{\|\bar L\|_1}\right)
\bigl(f(v_k)-f^\star\bigr)
\le
\left(1-\frac{\mu}{\|\bar L\|_1}\right)
\bigl(f(x_k)-f^\star\bigr).
\label{eq:inertial-linear}
\end{equation}
Consequently,
\begin{equation}
f(x_k)-f^\star
\le
\left(1-\frac{\mu}{\|\bar L\|_1}\right)^{k-1}
\bigl(f(x_1)-f^\star\bigr).
\label{eq:inertial-global}
\end{equation}
More precisely, whenever $f(x_k)>f^\star$, define
\[
\theta_k:=\frac{f(v_k)-f^\star}{f(x_k)-f^\star}\in[0,1].
\]
Then the one-step factor is bounded by
\begin{equation}
f(x_{k+1})-f^\star
\le
\left(1-\frac{\mu}{\|\bar L\|_1}\right)\theta_k
\bigl(f(x_k)-f^\star\bigr).
\label{eq:inertial-theta}
\end{equation}
Thus the restart gives the same worst-case linear guarantee as the basic adaptive SignGD step, while a successful extrapolation with $\theta_k<1$ improves that particular step. No uniform accelerated rate follows unless one can control $\theta_k$ away from~1.
\end{theorem}
\begin{proof}
Apply Proposition~\ref{prop:suff-decrease} at the extrapolated point $v_k$:
\[
f(x_{k+1})-f^\star
\le f(v_k)-f^\star-
\frac{\|\nabla f(v_k)\|_1^2}{2\|\bar L\|_1}.
\]
Strong convexity gives
$\|\nabla f(v_k)\|_1^2\ge2\mu(f(v_k)-f^\star)$, hence the first inequality in~\eqref{eq:inertial-linear}. By construction of the restart, $f(v_k)\le f(x_k)$, which gives the second. Iteration proves~\eqref{eq:inertial-global}; dividing the first inequality by the current gap gives~\eqref{eq:inertial-theta}.
\end{proof}

The theorem shows that the momentum/restart variant is globally safe under the same assumptions as the basic adaptive method and quantifies the gain of a successful extrapolation through $\theta_k$. Establishing a uniform factor $\theta_k<1$, or another accelerated complexity bound, requires additional analysis.

We now evaluate these methods on controlled convex benchmarks, comparing step policies and the inertial variant with/without restart.

\section{Numerical Experiments}
\label{sec:experiments}
The experiments have two goals. First, we check the descent behavior predicted by the certified $\ell_\infty$ analysis. Second, we test the inertial scheme with restart and the face-aware rules motivated by the Filippov discussion. We use one reproducible Python notebook for all experiments and export every paper figure as a vector PDF. The notebook is available at
\href{https://colab.research.google.com/drive/1q-QgXRuF7mHpA4wLST0DD_XsT2ntkgXi?usp=sharing}{\texttt{this Colab link}}.

\subsection{Experimental Protocol}
We consider three smooth strongly convex synthetic objectives. Their Hessians are uniformly bounded by a positive semidefinite matrix $G$ as summarized in Table~\ref{tab:benchmarks}. From $G$ we form the diagonal majorizer
\[
D=\operatorname{Diag}(L_1,\ldots,L_d),
\qquad
L_i=\sum_{j=1}^d |G_{ij}|,
\]
for which $D\succeq G$. This gives the valid separable model used by the face-restricted experiment. For the full SignGD step, we use the certified global bound
\begin{equation}
L_\infty
=\min\left\{\sum_{i=1}^d L_i,\;d\,\lambda_{\max}(G)\right\},
\label{eq:num-linf-certificate}
\end{equation}
which satisfies
$h^\top \nabla^2 f(x)h\le L_\infty\|h\|_\infty^2$ for all $x,h$.
This point is important: diagonal Hessian entries alone do not give the upper model required by the theory on dense problems.

\begin{table}[ht!]
\centering
\caption{Synthetic benchmark families and certified Hessian upper bounds}
\label{tab:benchmarks}
\resizebox{\textwidth}{!}{%
\begin{tabular}{@{}lll@{}}
\toprule
Problem & Objective & Hessian upper bound $G$ \\
\midrule
Logistic--Quadratic & $\frac12\|Ax\|_2^2+\gamma\sum_j\log(1+e^{(Bx)_j})$ & $A^\top A+\frac{\gamma}{4}B^\top B$ \\
Smooth Max & $\frac12x^\top Qx+\gamma\log\!\left(\sum_i e^{x_i}\right)$ & $Q+\frac{\gamma}{2}I$ \\
$\ell_2$-regularized logistic & $\frac1n\sum_j\log(1+e^{-y_j a_j^\top x})+\frac{\lambda}{2}\|x\|_2^2$ & $\frac{1}{4n}A^\top A+\lambda I$ \\
\bottomrule
\end{tabular}%
}
\end{table}

Unless stated otherwise, we use $d=200$, a budget of $2000$ iterations, and $10$ independent evaluation seeds. Each curve is the median across these runs and the shaded region is the interquartile range. Reference solutions are computed with L-BFGS to high accuracy. Values below roughly $10^{-12}$ should therefore be read as reference/numerical precision rather than as meaningful additional digits.

We compare five methods. The theory-certified SignGD step is
\[
\eta_k=\frac{\|\nabla f(x_k)\|_1}{L_\infty}.
\]
The inertial method uses the same step at its extrapolated point, a fixed momentum parameter $\beta=0.9$, and the monotone restart rule of Section~\ref{sec:inertial-signgd}. We also show the same inertial method without restart. Finally, we include constant-step SignGD and constant-step inertial SignGD. On the synthetic problems, the constant-step multiplier is selected on the disjoint calibration seeds $101,102,103$ and then frozen before evaluation on seeds $0,\ldots,9$. Hence the evaluation curves do not reuse the calibration instances. The real-data experiment uses a separate validation split for this selection.

\subsection{Main Synthetic Results}
\begin{figure}[ht!]
  \centering
  \includegraphics[width=0.48\linewidth]{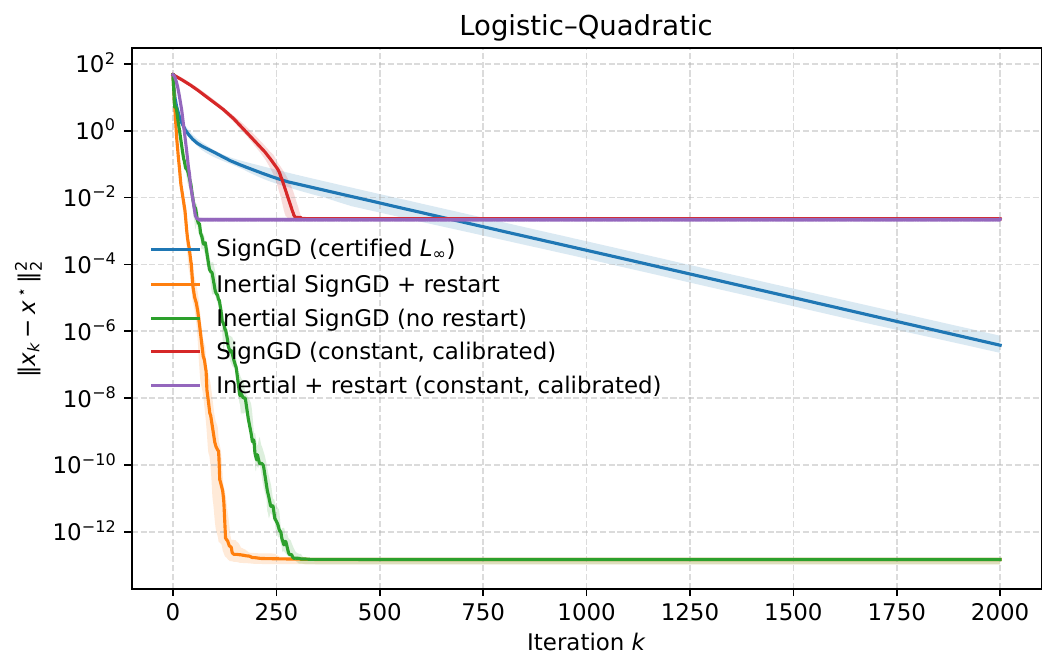}\hfill
  \includegraphics[width=0.48\linewidth]{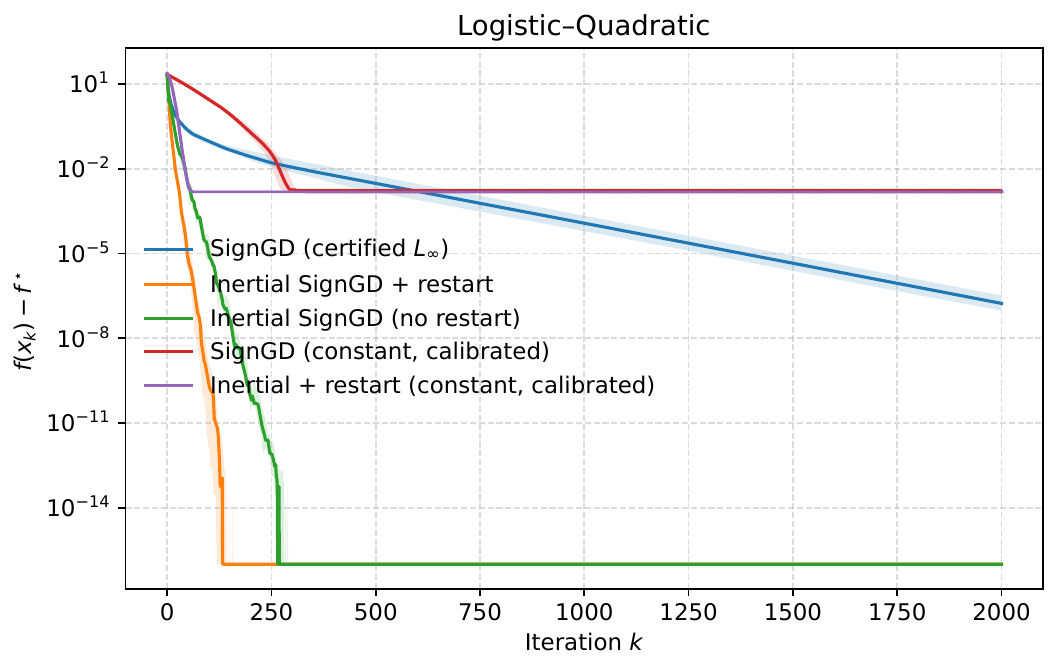}
  \caption{Logistic--Quadratic, $d=200$. Left: $\|x_k-x^\star\|_2^2$. Right: $f(x_k)-f^\star$. Curves show the median over $10$ evaluation seeds and shaded regions show the interquartile range. The adaptive methods use the certified $L_\infty$ in~\eqref{eq:num-linf-certificate}; constant-step multipliers are calibrated on disjoint seeds}
  \label{fig:lq}
\end{figure}

\begin{figure}[ht!]
  \centering
  \includegraphics[width=0.48\linewidth]{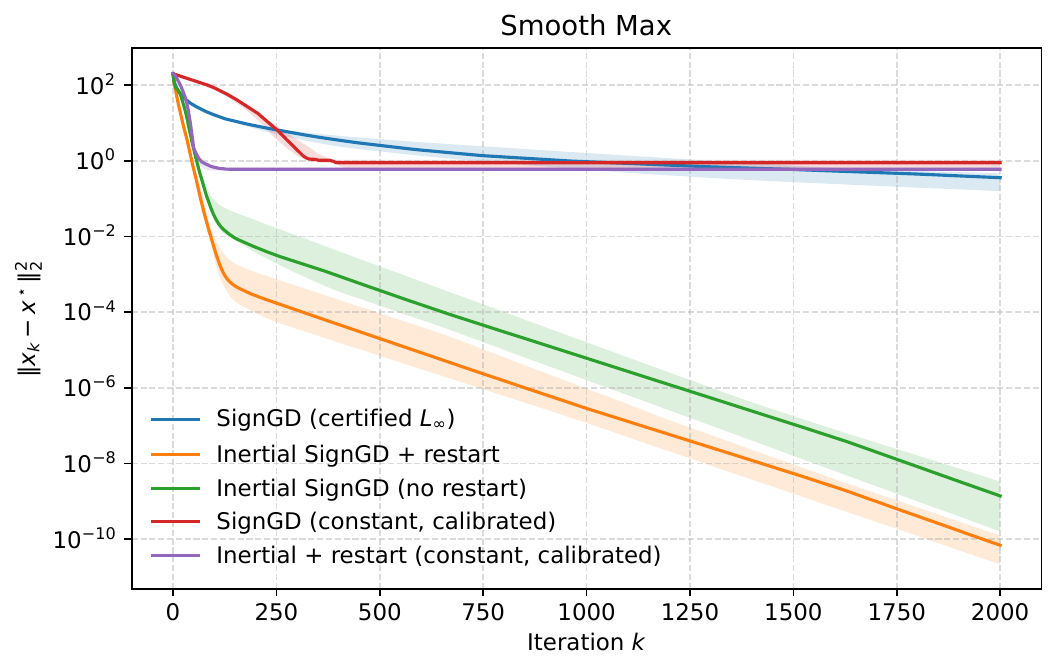}\hfill
  \includegraphics[width=0.48\linewidth]{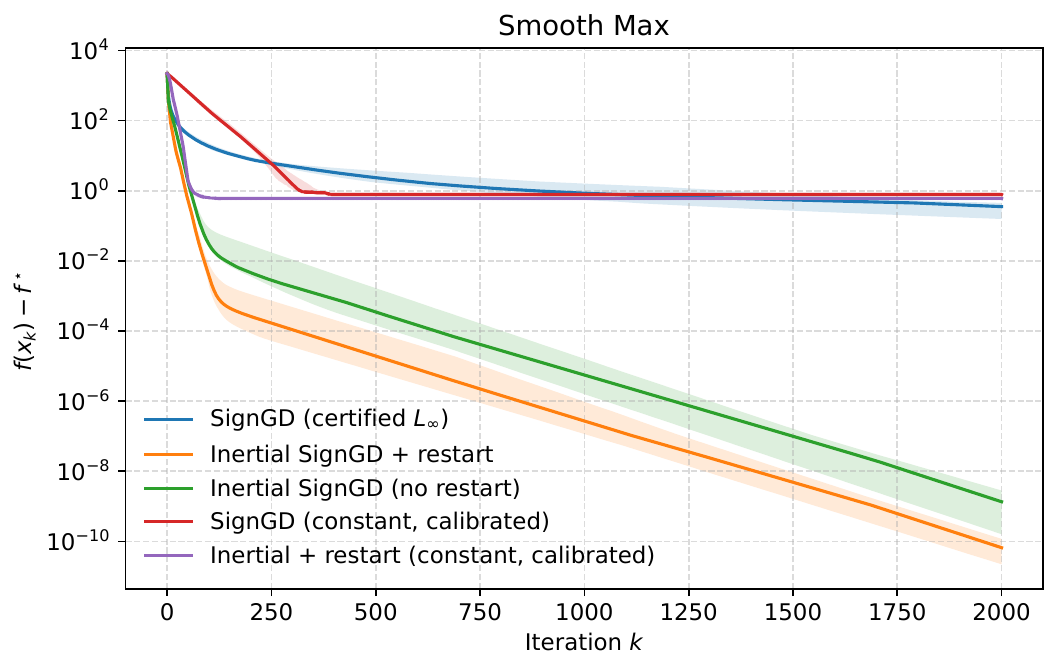}
  \caption{Smooth Max, $d=200$. Left: $\|x_k-x^\star\|_2^2$. Right: $f(x_k)-f^\star$. Median and interquartile range over $10$ evaluation seeds. The restarted inertial rule gives a large empirical speedup while using the same certified curvature bound}
  \label{fig:smax}
\end{figure}

\begin{figure}[ht!]
  \centering
  \includegraphics[width=0.48\linewidth]{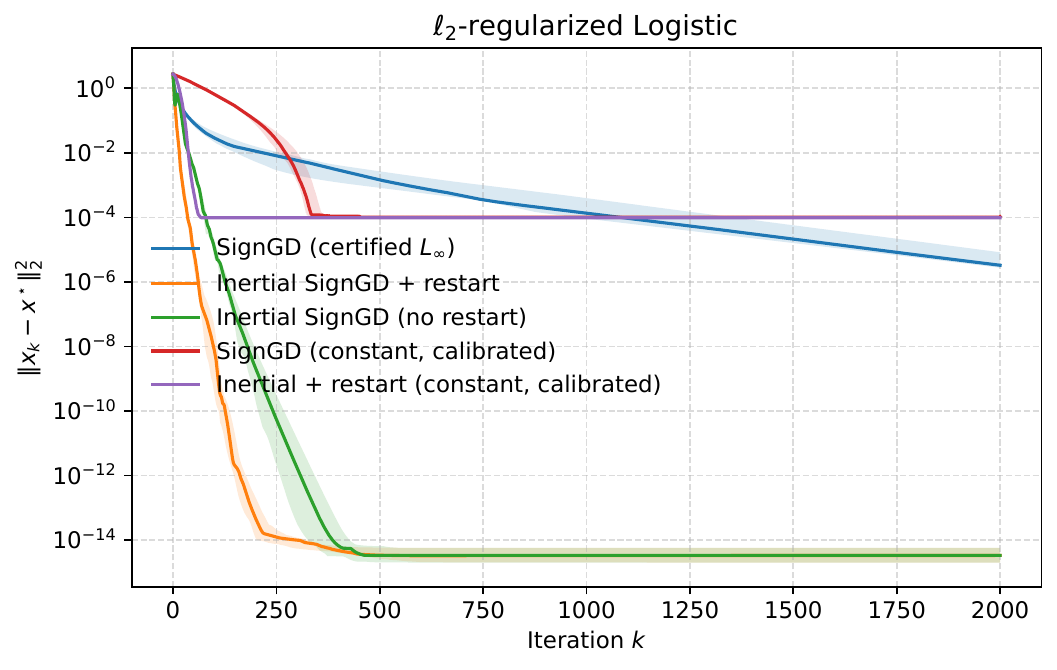}\hfill
  \includegraphics[width=0.48\linewidth]{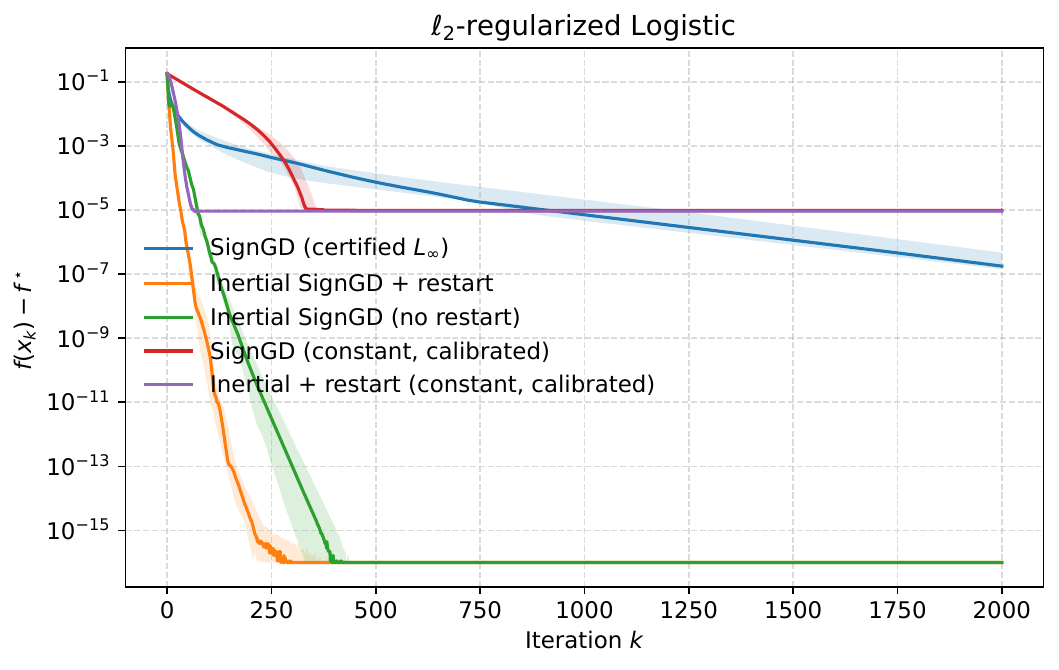}
  \caption{$\ell_2$-regularized logistic regression, $d=200$. Left: $\|x_k-x^\star\|_2^2$. Right: $f(x_k)-f^\star$. Median and interquartile range over $10$ evaluation seeds}
  \label{fig:logreg}
\end{figure}

Figures~\ref{fig:lq}--\ref{fig:logreg} show the same general behavior. The certified SignGD curves decrease monotonically on every run, as predicted by the descent estimate. The global constant $L_\infty$ is safe but can be conservative, so plain SignGD is sometimes slow. With the same certified curvature bound, the restarted inertial method reaches a small neighborhood of the reference solution much earlier on all three families.

Table~\ref{tab:num-summary} quantifies this effect. On Logistic--Quadratic and logistic regression, the restarted inertial method reaches a $10^{-4}$ function gap in about $41$ and $21$ median iterations, respectively, versus about $1026$ and $453$ for certified SignGD. On Smooth Max, certified SignGD does not reach $10^{-2}$ within the $2000$-iteration budget, whereas the restarted inertial method reaches $10^{-4}$ in about $299$ iterations. We emphasize that this is an empirical speedup; the theorem in Section~\ref{sec:inertial-signgd} gives a global linear guarantee but does not claim an accelerated worst-case complexity.

\begin{table}[ht!]
\centering
\caption{Median behavior over $10$ evaluation seeds. ``$>2000$'' means that the target was not reached within the budget. Gaps at reference precision are reported as $<10^{-12}$}
\label{tab:num-summary}
\small
\resizebox{\textwidth}{!}{%
\begin{tabular}{@{}lrrrr@{}}
\toprule
& \multicolumn{2}{c}{Median iterations to $f(x_k)-f^\star\le10^{-4}$} & \multicolumn{2}{c}{Median final gap} \\
Problem & SignGD & Inertial + restart & SignGD & Inertial + restart \\
\midrule
Logistic--Quadratic & $1026$ & $41$ & $1.7\times10^{-7}$ & $<10^{-12}$ \\
Smooth Max & $>2000$ & $299$ & $3.5\times10^{-1}$ & $6.6\times10^{-11}$ \\
$\ell_2$-regularized logistic & $453$ & $21$ & $1.8\times10^{-7}$ & $<10^{-12}$ \\
\bottomrule
\end{tabular}%
}
\end{table}

The restart is also useful. On iterations for which the Logistic--Quadratic gap is above $10^{-10}$, successful extrapolations have median $\theta_k\simeq0.88$, where
\[
\theta_k=\frac{f(v_k)-f^\star}{f(x_k)-f^\star}.
\]
About $8\%$ of the attempted extrapolations are rejected. This directly illustrates the factor in the inertial convergence estimate: accepted extrapolations lower the starting gap of the next SignGD step, while restart removes unfavorable ones. The corresponding diagnostic is reported in Appendix~\ref{app:ablation_active_face}.

\paragraph{Logistic--Quadratic.}
We use $n_A=n_B=800$, $d=200$, $\gamma=1$, column-normalized Gaussian matrices $A,B$, and $x_0=0$. The restarted inertial method reaches reference precision rapidly. The no-restart variant also converges, but it needs roughly twice as many iterations to reach the smallest displayed tolerances. The calibrated constant steps decrease quickly at first and then settle at the expected fixed-step neighborhood, while the gradient-scaled rules continue to reduce the gap.

\paragraph{Smooth Max.}
We use $d=200$, $\gamma=1$, and $Q=U\operatorname{Diag}(\lambda)U^\top$ with eigenvalues logarithmically spaced between $1$ and $100$. This is the most conservative case for the global $L_\infty$ bound. Plain certified SignGD therefore makes steady but slow progress. In contrast, the inertial variants preserve the same certified local step and reduce the gap by many additional orders of magnitude. Restart gives a clear improvement over the no-restart trajectory.

\paragraph{$\ell_2$-regularized logistic regression.}
We generate $n=2000$ samples in dimension $d=200$, standardize the features, and use $\lambda=10^{-3}$. Again, the restarted inertial rule reaches reference precision much sooner than certified SignGD. The constant-step methods reach a visible error floor, whereas the scaled step continues to shrink with the gradient.

\subsection{Numerical Evidence for the Face-Aware Results}
We next examine the face-aware results of Sections~\ref{sec:signgd-strong} and~\ref{sec:filippov} directly.

First, the mass-aware restriction from Proposition~\ref{prop:mass-aware-face} is tested on the Logistic--Quadratic family. Keeping only enough coordinates to retain $90\%$ of the $\ell_1$ gradient mass uses a median of about $50$ coordinates out of $200$ and only about $25\%$ of the full separable curvature sum. Its median final gap is $9.9\times10^{-8}$, compared with $2.5\times10^{-3}$ for the full separable-majorizer step. Retaining $99\%$ of the gradient mass uses about $152$ coordinates and $76\%$ of the curvature sum, and gives the intermediate final gap $7.5\times10^{-4}$. This illustrates the tradeoff between retained gradient mass and active curvature in Proposition~\ref{prop:mass-aware-face}. The complete curves and diagnostics are given in Appendix~\ref{app:ablation_active_face}.

Second, we compare ordinary SignGD with the one-hit and two-hit face-aware rules on synthetic logistic regression in dimensions $20$ and $100$. In these tests, a strict sign crossing requires two consecutive partial derivatives to have opposite signs and magnitude greater than $10^{-14}$. All variants converge to the same reference solution. One-hit is slightly faster on these regular instances; the two-hit rule is slightly slower but remains very close to SignGD. More importantly for Theorem~\ref{thm:safeguarded-sliding}, with $\rho=0.5$ the safeguard accepts every two-hit candidate in all $20$ evaluation runs. At $d=100$, the median candidate score on event iterations is about $0.90$ and the smallest observed score is about $0.78$, well above the threshold. The largest residual in the certified descent inequality is at roundoff level ($2.2\times10^{-16}$). Thus the safeguard certifies the raw two-hit steps in these experiments without changing their trajectory. Full results are reported in Appendix~\ref{app:test_Proj_SGD}.

\subsection{Real-Data Check}
Finally, Appendix~\ref{app:realdata} reports the Breast Cancer Wisconsin diagnostic data set. The aim is only to check that the same behavior persists on a small real classification problem. The certified inertial rule reduces the training objective much more rapidly than certified SignGD, while validation-tuned constant steps give slightly better held-out loss. All tested methods obtain test accuracy between $98.2\%$ and $99.1\%$. This distinction is expected: our theory concerns optimization, not statistical generalization.

\subsection{Summary of the Numerical Evidence}
The experiments support four points used in the paper.
\begin{itemize}
\item The certified SignGD step is monotone on all $30$ synthetic evaluation runs (three problem families and ten seeds each).
\item Inertial SignGD with restart gives a large and consistent empirical reduction in iteration counts while keeping the same certified SignGD step. The restart diagnostic agrees with the factor $\theta_k$ appearing in the convergence bound.
\item The mass-aware experiment illustrates the practical meaning of the face-restricted rate: retaining most of the gradient mass while discarding a large part of the curvature can substantially improve the step.
\item The one-hit and two-hit rules remain stable under frequent switching. For the two-hit rule, every candidate in these tests satisfies the $\rho=0.5$ safeguard, and every numerical certificate is satisfied to machine precision.
\end{itemize}
The notebook also contains a direct sufficient-decrease backtracking experiment. It improves plain SignGD over the conservative global $L_\infty$ bound on Logistic--Quadratic, as expected. We report it only as a reproducibility diagnostic, not as a separate algorithmic contribution.

\section{Conclusions}

We studied SignGD through norm-constrained continuous-time dynamics. The optimal-velocity formulation makes zeros and ties explicit by keeping the full set of norm-steepest directions. This leads to existence of trajectories, an exact descent identity, convergence estimates for convex objectives, and finite-time convergence under a Polyak--\L{}ojasiewicz inequality expressed in the dual norm.

Filippov regularization describes the dynamics obtained after selecting a single-valued direction from the optimal-velocity set. Away from switching sets, the selector follows an ordinary differential equation. At a switch, the Filippov set collects the locally attainable convexified velocities. This gives a simple local criterion for crossing and sliding and motivates the one-hit freeze, two-hit sliding-track, and tie-aware convex-combination updates.

The discrete analysis gives descent certificates for damped sign directions and a safeguard that preserves a global linear rate under the corresponding smoothness and Polyak--\L{}ojasiewicz assumptions. The same framework gives a linear guarantee for tie-aware $\ell_1$ sliding steps. Under separable smoothness and strong convexity, the mass-aware analysis further shows how a face can be selected by balancing retained gradient mass against the curvature carried by the selected coordinates. The restarted inertial scheme preserves the global linear guarantee of the adaptive SignGD step and can improve individual iterations when extrapolation is successful.

The numerical experiments support these conclusions. The certified SignGD step is monotone on all synthetic evaluation runs. The restarted inertial method substantially reduces iteration counts on the tested problems, and the mass-aware experiment illustrates the predicted gradient-mass--curvature tradeoff. The two-hit safeguard is satisfied throughout the reported logistic-regression tests. Natural directions for further study include sharper guarantees for the inertial scheme, stochastic versions of the face-aware dynamics, and a more detailed comparison between velocity-constrained flows and norm-ball linear-minimization-oracle methods.

\appendix
\section{Additional Details}


\subsection{Steepest Descent Under Norm Constraints (Details)}
\label{app:dual-norm-lemma}

\paragraph{Lemma A.1 (Steepest descent under a norm via the dual norm).}
Let $\|\cdot\|$ be a norm on $\R^d$ with dual norm
$\|y\|_\ast := \sup_{\|\vv\|\le 1} \langle y,\vv\rangle$.
For any $g\in\R^d$,
\begin{equation}\label{eq:app-min-equals-minus-dual}
\min_{\|\vv\|\le 1}\ \langle g,\vv\rangle \;=\; -\,\|g\|_\ast.
\end{equation}
Moreover, the set of minimizers is
\begin{equation}\label{eq:app-argmin-is-neg-subgrad}
\arg\min_{\|\vv\|\le 1}\ \langle g,\vv\rangle \;=\; -\,\partial\|\cdot\|_\ast(g),
\end{equation}
where $\partial\|\cdot\|_\ast(g)$ is the subdifferential of the convex function
$y\mapsto \|y\|_\ast$ at $g$.

\emph{Proof.}
By definition of the dual norm,
$\|g\|_\ast = \sup_{\|\vv\|\le 1}\langle g,\vv\rangle$.
Then
\[
\min_{\|\vv\|\le 1}\langle g,\vv\rangle
= -\,\sup_{\|\vv\|\le 1}\langle -g,\vv\rangle
= -\,\|-g\|_\ast
= -\,\|g\|_\ast,
\]
using symmetry of norms. For the argmin, recall the subgradient characterization of a norm:
\begin{equation}\label{eq:app-norm-subgrad-char}
s\in \partial \|{\cdot}\|_\ast(g)
\quad\Longleftrightarrow\quad
\|s\|\le 1\ \ \text{and}\ \ \langle s,g\rangle=\|g\|_\ast.
\end{equation}
A feasible $v^\star$ with $\|v^\star\|\le 1$ achieves \eqref{eq:app-min-equals-minus-dual} iff
$\langle g,v^\star\rangle=-\|g\|_\ast$.
Setting $s^\star:=-v^\star$ gives $\|s^\star\|\le 1$ and $\langle s^\star,g\rangle=\|g\|_\ast$,
so $s^\star\in\partial\|\cdot\|_\ast(g)$ by \eqref{eq:app-norm-subgrad-char}.
Equivalently, $v^\star=-s^\star$ with $s^\star\in\partial\|\cdot\|_\ast(g)$, proving
\eqref{eq:app-argmin-is-neg-subgrad}. \hfill$\square$

\paragraph{Concrete subgradients (worked out).}
\begin{itemize}
\item \textbf{$\ell_2$ case.} For $g\neq 0$, $\partial\|\cdot\|_2(g)=\{g/\|g\|_2\}$, hence
$\arg\min=\{-g/\|g\|_2\}$ (normalized GD). If $g=0$, any $\|\vv\|_2\le 1$ minimizes.

\item \textbf{$\ell_\infty$ constraint (dual $\ell_1$).}
\[
\partial\|\cdot\|_1(g)=\Big\{s\in\R^d:\ s_i=\operatorname{sign}(g_i)\ \text{if }g_i\neq 0,\ 
\ s_i\in[-1,1]\ \text{if }g_i=0\Big\}.
\]
Thus $\arg\min_{\|\vv\|_\infty\le 1}\langle g,\vv\rangle=-\partial\|\cdot\|_1(g)$:
componentwise $v_i^\star=-\operatorname{sign}(g_i)$, with $v_i^\star\in[-1,1]$ when $g_i=0$.

\item \textbf{$\ell_1$ constraint (dual $\ell_\infty$).}
If $g\neq0$, let $I(g):=\arg\max_i|g_i|$. Then
\[
\partial\|\cdot\|_\infty(g)
=\operatorname{conv}\Big\{\operatorname{sign}(g_i)\,\ve^{(i)}:i\in I(g)\Big\}.
\]
Thus any convex combination of the tied signed basis directions is optimal; choosing one extreme point gives the usual greedy coordinate step. If $g=0$, then
$\partial\|\cdot\|_\infty(0)=\{v:\|v\|_1\le1\}$, so the whole $\ell_1$ unit ball is optimal.
\end{itemize}

\paragraph{Geometric picture (support functions).}
The dual norm $\|g\|_\ast$ is the support function of the primal unit ball
$B:=\{\vv:\|\vv\|\le 1\}$. Minimizers of $\langle g,\cdot\rangle$ over $B$ form the face of $B$ on which this linear form is smallest: a vertex (unique direction) or a higher-dimensional face (a convex set of directions). The latter case corresponds to ties/zeros and leads to set-valued dynamics---handled rigorously via Filippov convexification in Section~\ref{sec:filippov}.

\subsection{Mass-Aware Face Restriction and Inertial Diagnostic}\label{app:ablation_active_face}
This appendix gives a direct numerical test of Proposition~\ref{prop:mass-aware-face}. The aim is to measure the tradeoff between retained gradient mass and the curvature paid by the selected coordinates.

For the Logistic--Quadratic family, let $J_k$ be the smallest set of coordinates, ordered by decreasing $|g_{k,i}|$, such that
\[
\theta_k=\frac{\sum_{i\in J_k}|g_{k,i}|}{\|g_k\|_1}\ge\theta_0.
\]
We then use the certified face-restricted step
\[
A_k=\sum_{i\in J_k}|g_{k,i}|,
\qquad
S(J_k)=\sum_{i\in J_k}L_i,
\qquad
x_{k+1}=x_k-\frac{A_k}{S(J_k)}d_k,
\]
where $(d_k)_i=\operatorname{sign}(g_{k,i})$ on $J_k$ and zero otherwise. We compare $\theta_0\in\{0.90,0.99\}$ with the full separable-majorizer step $J_k=\{1,\ldots,d\}$. All curves aggregate the same $10$ evaluation seeds as Section~\ref{sec:experiments}.

\begin{figure}[ht!]
  \centering
  \includegraphics[width=0.78\linewidth]{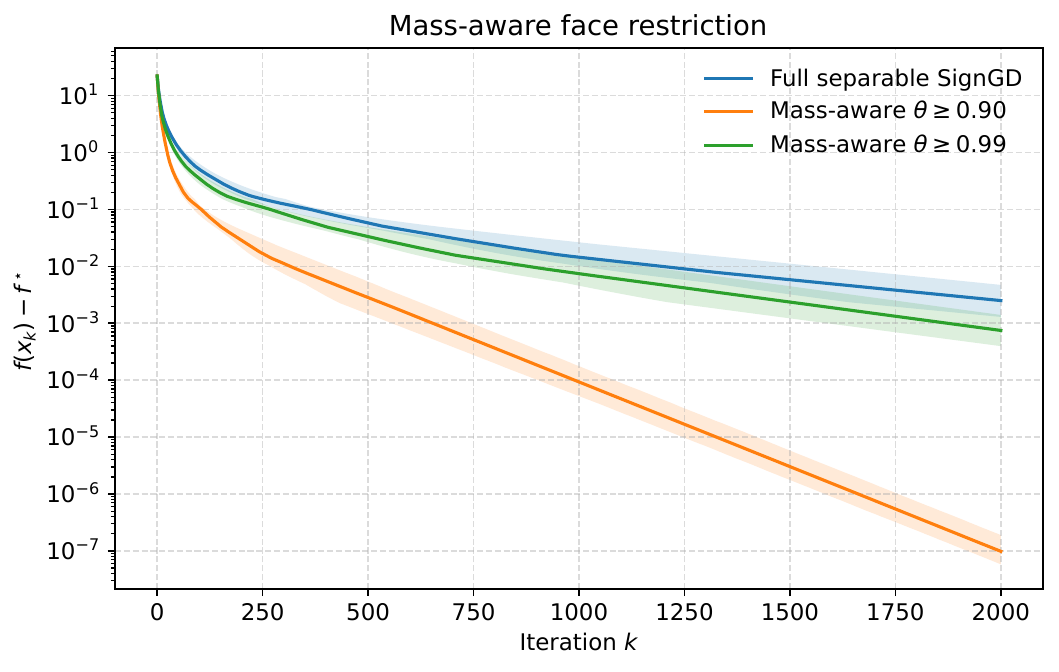}
  \caption{Mass-aware face restriction on Logistic--Quadratic. Curves are medians over $10$ seeds and bands show the interquartile range. Retaining $90\%$ of the gradient $\ell_1$ mass gives the strongest reduction because it removes a much larger fraction of the separable curvature}
  \label{fig:mass-aware-gap}
\end{figure}

Figure~\ref{fig:mass-aware-gap} shows a clear curvature--mass tradeoff. With $\theta_0=0.90$, the method keeps a median of about $50$ coordinates out of $200$ and only $25\%$ of the full curvature sum, while retaining just over $90\%$ of the gradient $\ell_1$ mass. The median final gap is $9.9\times10^{-8}$, versus $2.5\times10^{-3}$ for the full separable step. The more conservative choice $\theta_0=0.99$ keeps about $152$ coordinates and $76\%$ of the curvature sum and ends at $7.5\times10^{-4}$. These results agree with Proposition~\ref{prop:mass-aware-face}: the useful quantity is not support size alone, but the amount of descent mass retained relative to the curvature paid on that support.

\begin{figure}[ht!]
  \centering
  \includegraphics[width=\linewidth]{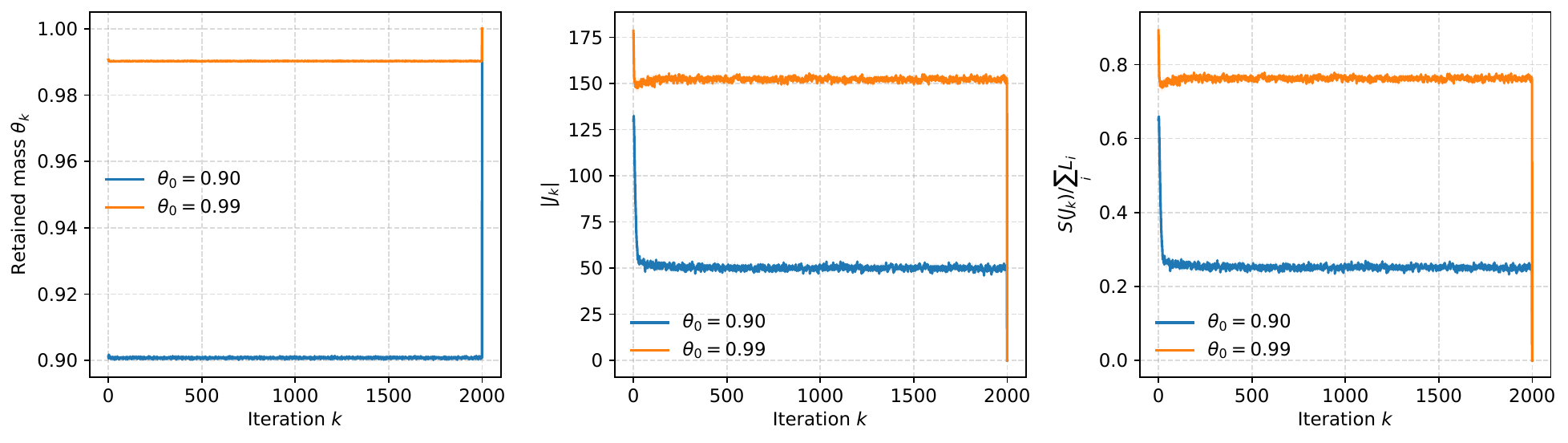}
  \caption{Diagnostics for the mass-aware experiment. Left: retained gradient mass $\theta_k$. Middle: selected support size $|J_k|$. Right: curvature ratio $S(J_k)/\sum_iL_i$. For $\theta_0=0.90$, about one quarter of the total separable curvature is used while about $90\%$ of the gradient mass is retained}
  \label{fig:mass-aware-diagnostics}
\end{figure}

The experiment also shows why the mass-aware statement is useful on dense problems. Exact zero-gradient coordinates are not needed; what matters is the fraction of gradient mass kept by the selected face.

\begin{figure}[ht!]
  \centering
  \includegraphics[width=\linewidth]{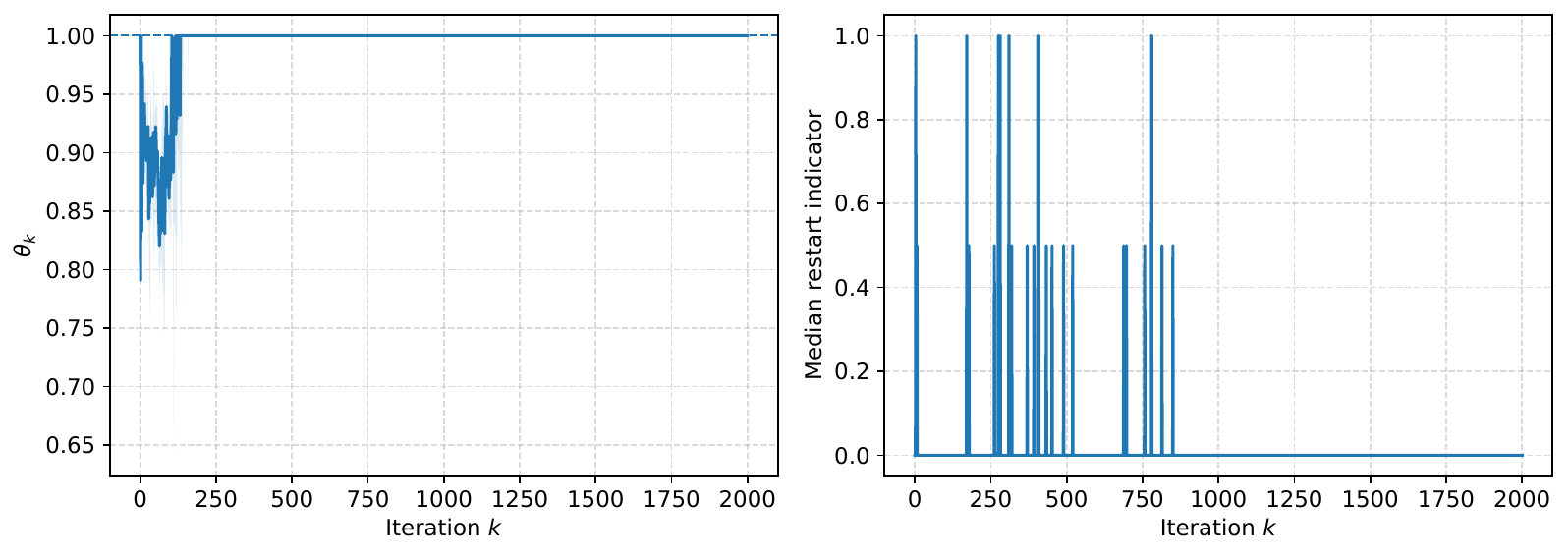}
  \caption{Inertial SignGD with restart on Logistic--Quadratic. Left: median $\theta_k$ and interquartile range. Right: median restart indicator. During the pre-asymptotic phase, successful extrapolations typically satisfy $\theta_k<1$, while restart removes extrapolations that increase the objective}
  \label{fig:inertial-theta}
\end{figure}

\paragraph{Inertial restart diagnostic.}
Figure~\ref{fig:inertial-theta} reports the quantity
\[
\theta_k=\frac{f(v_k)-f^\star}{f(x_k)-f^\star}
\]
This quantity appears in the inertial convergence bound. Figure~\ref{fig:inertial-theta} also shows the restart indicator for the Logistic--Quadratic test. While the gap is above $10^{-10}$, successful extrapolations have median $\theta_k\simeq0.88$, with an interquartile range of approximately $0.83$--$0.92$. About $8\%$ of attempted extrapolations are rejected. Once numerical precision is reached, the displayed ratio is set to $1$ by convention. Thus the figure describes the pre-asymptotic phase.

\subsection{Face-Aware SignGD: One-Hit and Two-Hit Rules}
\label{app:test_Proj_SGD}
We compare ordinary scaled SignGD with the two face-aware rules introduced in Section~\ref{sec:filippov}: one-hit freeze and two-hit sliding-track. For the two-hit rule we also report the safeguarded version of Theorem~\ref{thm:safeguarded-sliding}, with the fixed threshold $\rho=0.5$.

The objective is $\ell_2$-regularized logistic regression,
\[
\ell(w)=\frac1n\sum_{i=1}^n\log\bigl(1+e^{-y_iw^\top x_i}\bigr)+\frac{\lambda}{2}\|w\|_2^2,
\]
with $n=2000$, $\lambda=10^{-3}$, and dimensions $d\in\{20,100\}$. Features are standardized. Every method uses the same certified $L_\infty$ step, and the plots aggregate $10$ independent seeds over $2000$ iterations. A sign change is counted only when both consecutive partial derivatives have magnitude greater than $10^{-14}$ and have opposite signs. The same tolerance is used in the two-hit interpolation: when $|\beta|\le10^{-14}$, the local model is treated as degenerate and the default sign direction is used. The two-hit construction is exactly the one in Algorithm~\ref{alg:two_hit_sliding}; the safeguard changes the candidate only when its score
\[
q(w,g)=\frac{2\sum_iw_i|g_i|}{\|g\|_1}-\|w\|_\infty^2
\]
falls below $\rho$.

\begin{figure}[ht!]
  \centering
  \includegraphics[width=0.78\textwidth]{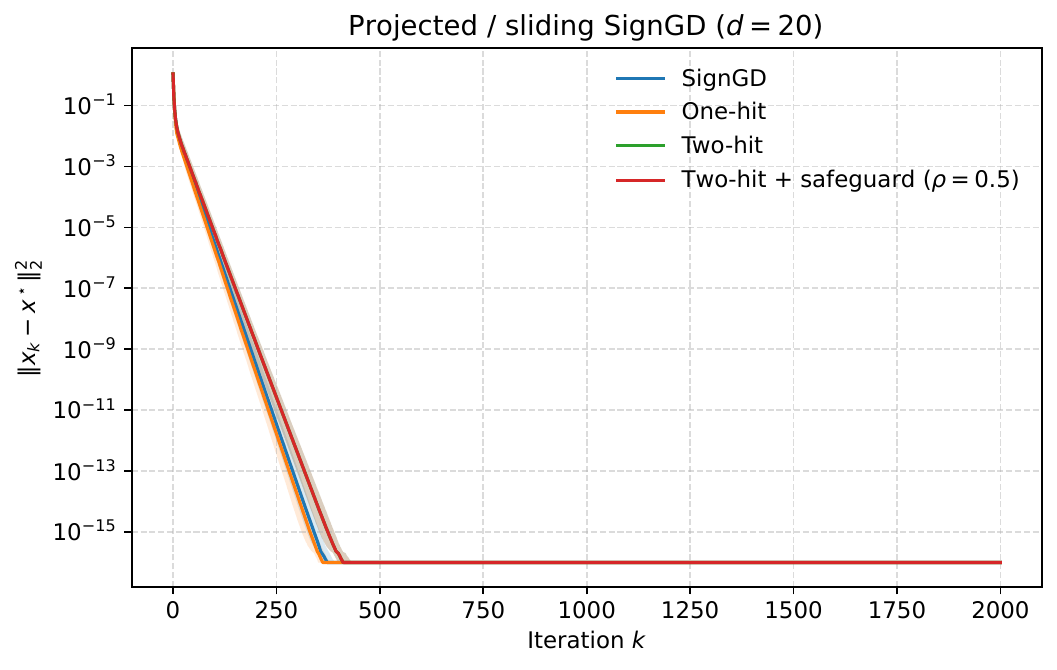}
  \caption{Squared distance $\|w_k-w^\star\|_2^2$ for $d=20$. Curves show medians over $10$ seeds. One-hit gives a small improvement; raw and safeguarded two-hit coincide}
  \label{fig:d20}
\end{figure}

\begin{figure}[ht!]
  \centering
  \includegraphics[width=0.78\textwidth]{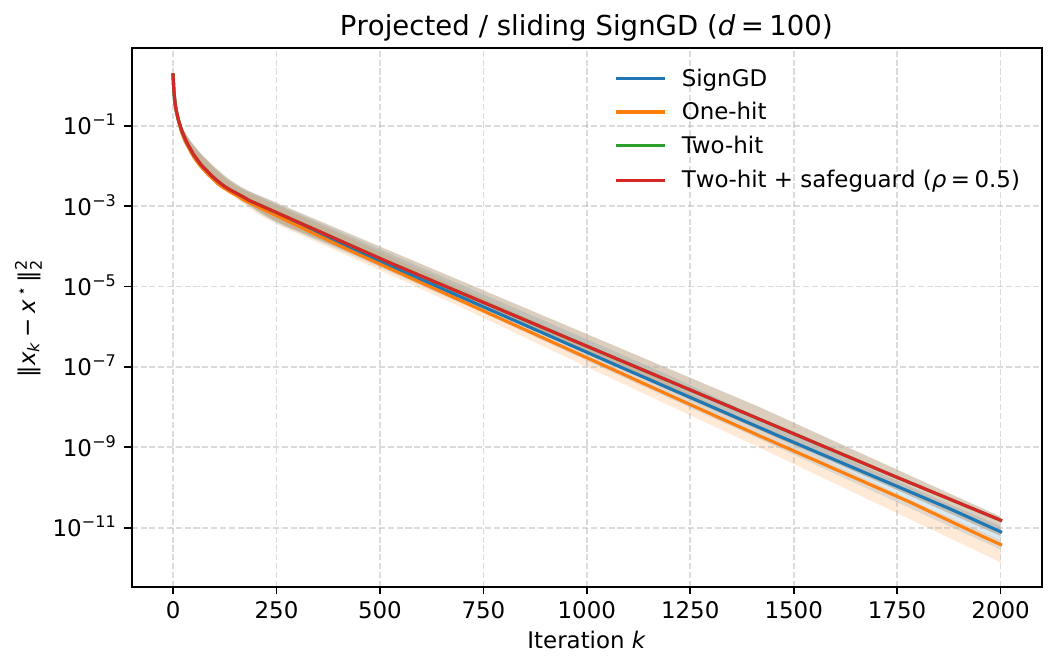}
  \caption{Squared distance $\|w_k-w^\star\|_2^2$ for $d=100$. The four methods remain close over the full run; one-hit is slightly faster, while the safeguard leaves the two-hit trajectory unchanged}
  \label{fig:d100}
\end{figure}

\paragraph{Convergence behavior.}
Figures~\ref{fig:d20} and~\ref{fig:d100} show that all four trajectories converge to the same reference solution. The differences are small but systematic. One-hit is slightly faster than ordinary SignGD, while two-hit is slightly slower. For example, at $d=100$ the median iteration counts to reach a $10^{-8}$ function gap are approximately $1013$ (one-hit), $1043$ (SignGD), and $1080$ (two-hit). The raw and safeguarded two-hit curves coincide in both dimensions.

On these regular logistic instances, the safeguard does not reject any two-hit candidates. Its role here is therefore only to certify the raw two-hit steps. Across all $20$ runs, no two-hit candidate is rejected. At $d=20$ the median candidate score on event iterations is $0.934$ and the smallest observed value is $0.598$; at $d=100$ the corresponding values are $0.901$ and $0.777$. All are above $\rho=0.5$. The maximum numerical residual in the sufficient-decrease inequality is $2.2\times10^{-16}$.

\begin{figure}[ht!]
  \centering
  \includegraphics[width=0.78\textwidth]{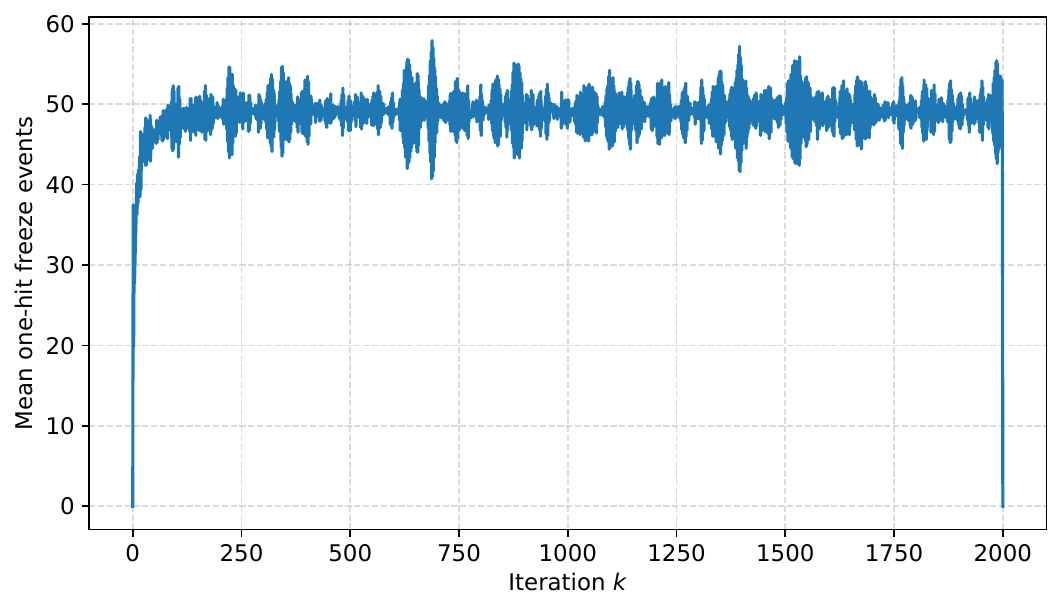}
  \caption{One-hit rule, $d=100$: mean number of frozen coordinates per iteration, averaged over $10$ seeds}
  \label{fig:flip_signs_1hit}
\end{figure}

\begin{figure}[ht!]
  \centering
  \includegraphics[width=0.78\textwidth]{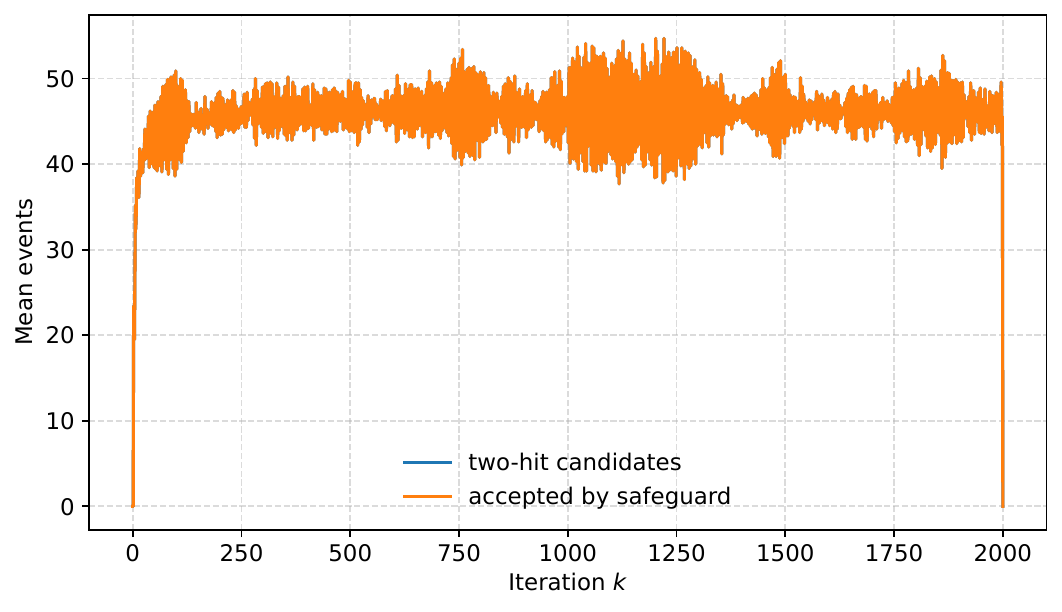}
  \caption{Two-hit rule, $d=100$: mean number of candidate coordinates and mean number accepted by the safeguard. The curves overlap: the safeguard does not reject any candidate in these runs}
  \label{fig:flip_signs_2hit}
\end{figure}

\paragraph{Switching statistics.}
At $d=100$, sign changes are frequent near the optimum: one-hit freezes a median of about $49$ coordinates per iteration over the complete run, while the two-hit detector produces about $46$ candidate coordinates per iteration. Hence the close convergence curves are not explained by an absence of switching. Figure~\ref{fig:flip_signs_1hit} shows the one-hit activity, while Figure~\ref{fig:flip_signs_2hit} compares two-hit candidates with candidates accepted by the safeguard. The two curves in the latter figure overlap because every candidate satisfies $q\ge\rho$.

These generic logistic instances do not show an empirical advantage for the more elaborate two-hit interpolation. They show, however, that the construction remains stable when switching is frequent and that the safeguard certifies the two-hit steps without modifying their trajectories.

\newpage
\subsection{Real-Data Logistic Regression}\label{app:realdata}
We use the Breast Cancer Wisconsin diagnostic data set distributed with \texttt{scikit-learn}. The split is fixed and stratified: $60\%$ training, $20\%$ validation, and $20\%$ test. Features are standardized using the training set and we use ridge parameter $\lambda=10^{-3}$. The optimization curves use the training objective. Constant steps are selected on the validation split using the same fixed tuning budget as in the notebook, and the test set is evaluated only after this selection.

\begin{figure}[ht!]
  \centering
  \includegraphics[width=\linewidth]{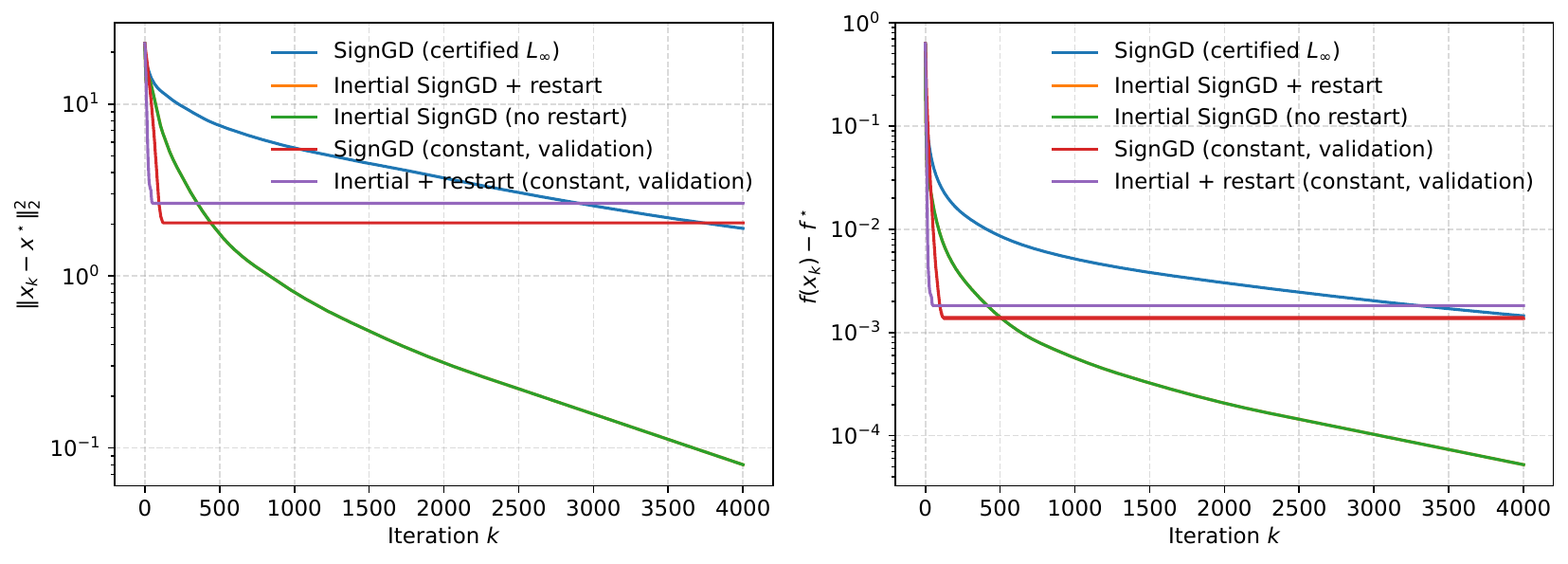}
  \caption{$\ell_2$-regularized logistic regression on the Breast Cancer data set. Left: $\|x_k-x^\star\|_2^2$. Right: $f(x_k)-f^\star$. Constant steps are selected using validation loss; the adaptive methods use the certified $L_\infty$ step}
  \label{fig:real-logreg}
\end{figure}

Figure~\ref{fig:real-logreg} shows the optimization trajectories. Certified inertial SignGD reduces the training gap more rapidly than certified SignGD. In this instance the restart condition is never triggered, so the restarted and no-restart adaptive inertial curves coincide. This is not a problem: restart is only a safeguard and need not be active on every problem.

Table~\ref{tab:real-logreg} reports optimization and prediction metrics separately. Among the adaptive methods, the inertial scheme has the smallest training gap. Validation-tuned constant SignGD obtains the smallest test loss. Test accuracies differ by only one sample on this small test set. We therefore do not interpret the experiment as evidence of statistical superiority; it shows only that the optimization behavior of the proposed schemes is not restricted to synthetic objectives.

\begin{table}[ht!]
\centering
\caption{Real-data logistic regression. The no-restart adaptive inertial run is omitted because it coincides with the restarted run on this data set}
\label{tab:real-logreg}
\small
\resizebox{\textwidth}{!}{%
\begin{tabular}{@{}lrrrr@{}}
\toprule
Method & Train gap & Validation loss & Test loss & Test accuracy \\
\midrule
SignGD, certified $L_\infty$ & $1.45\times10^{-3}$ & $0.0462$ & $0.0726$ & $98.25\%$ \\
Inertial SignGD + restart & $5.2\times10^{-5}$ & $0.0469$ & $0.0732$ & $99.12\%$ \\
SignGD, validation-tuned constant & $1.37\times10^{-3}$ & $0.0449$ & $\mathbf{0.0672}$ & $99.12\%$ \\
Inertial + restart, validation-tuned constant & $1.82\times10^{-3}$ & $\mathbf{0.0437}$ & $0.0713$ & $98.25\%$ \\
\bottomrule
\end{tabular}%
}
\end{table}

\bibliographystyle{spmpsci}
\bibliography{references}

\end{document}